\input amstex
\magnification=1200
\loadmsbm
\loadeufm
\loadeusm
\UseAMSsymbols
\input amssym.def

\font\BIGtitle=cmr10 scaled\magstep3
\font\bigtitle=cmr10 scaled\magstep1
\font\boldsectionfont=cmb10 scaled\magstep1
\font\section=cmsy10 scaled\magstep1

\def\scr#1{{\fam\eusmfam\relax#1}}

\def\scrD{{\scr D}}

\def\scrF{{\scr F}}

\def\scrL{{\scr L}}
\def\scrK{{\scr K}}

\def\scrO{{\scr O}}

\def\scrR{{\scr R}}
\def\scrT{{\scr T}}

\def\scrZ{{\scr Z}}
\def\gr#1{{\fam\eufmfam\relax#1}}

	\def\grg{{\gr g}}

	\def\grj{{\gr j}}
\def\grK{{\gr K}}	
	
\def\grl{{\gr l}}
	
	\def\grn{{\gr n}}

	\def\grs{{\gr s}}

	\def\grv{{\gr v}}

\def\db#1{{\fam\msbfam\relax#1}}

\def\dbC{{\db C}} 
 \def\dbF{{\db F}}
\def\dbG{{\db G}}

 \def\dbN{{\db N}}

 \def\dbZ{{\db Z}}

\def\gl{{\grg\grl}}

\def\eps{{\varepsilon}}

\def\char{\text{char}}
\def\Ker{\text{Ker}}
\def\n{\text{n}}

\def\Aut{\text{Aut}}
\def\der{\text{der}}

\def\ab{\text{ab}}

\def\ad{\text{ad}}
\def\Ad{\text{Ad}}

\def\Hom{\text{Hom}}
\def\End{\text{End}}
\def\Spec{\text{Spec}\,}

\def\sc{\text{sc}}
\def\Lie{\text{Lie}}

\def\leaderfill{\leaders\hbox to 1em
     {\hss.\hss}\hfill}
\def\nspace{\lineskip=1pt\baselineskip=12pt\lineskiplimit=0pt}

\def\finishproclaim{\par\rm
     \ifdim\lastskip<\medskipamount\removelastskip
     \penalty55\medskip\fi}
\def\endproof{$\hfill \square$}
\def\proof{\par\noindent {\it Proof:}\enspace}
\def\references#1{\par
  \centerline{\boldsectionfont References}\medskip
     \parindent=#1pt\nspace}
\def\Ref[#1]{\par\hang\indent\llap{\hbox to\parindent
     {[#1]\hfil\enspace}}\ignorespaces}
\def\Item#1{\par\smallskip\hang\indent\llap{\hbox to\parindent
     {#1\hfill$\,\,$}}\ignorespaces}
\def\ItemItem#1{\par\indent\hangindent2\parindent
     \hbox to \parindent{#1\hfill\enspace}\ignorespaces}

\def\Ge{{\mathchoice{\,{\scriptstyle\ge}\,}
  {\,{\scriptstyle\ge}\,}
  {\,{\scriptscriptstyle\ge}\,}{\,{\scriptscriptstyle\ge}\,}}}

\def\arrowsim{\,\smash{\mathop{\to}\limits^{\lower1.5pt
  \hbox{$\scriptstyle\sim$}}}\,}

\def\doublemaprights#1#2#3#4{\raise3pt\hbox{$\mathop{\,\,\hbox to     
#1pt{\rightarrowfill}\kern-30pt\lower3.95pt\hbox to
     #2pt{\rightarrowfill}\,\,}\limits_{#3}^{#4}$}}

\def\rightcapdownarrow{\raise9pt\hbox{$\ssize\cap$}\kern-7.75pt
     \Big\downarrow}

\def\rcapmapdown#1{\rightcapdownarrow\kern-1.0pt\vcenter{
     \hbox{$\scriptstyle#1$}}}

\def\rmapdown#1{\Big\downarrow\kern-1.0pt\vcenter{
     \hbox{$\scriptstyle#1$}}}
\def\rightsubsetarrow#1{{\ssize\subset}\kern-4.5pt\lower2.85pt
     \hbox to #1pt{\rightarrowfill}}
\def\longtwoheadedrightarrow#1{\raise2.2pt\hbox to #1pt{\hrulefill}
     \!\!\!\twoheadrightarrow}

\def\Hom{\operatorname{\hbox{Hom}}}

\def\im{\hbox{Im}}

\NoBlackBoxes
\parindent=25pt
\document
\footline={\hfil}

\null
\noindent 

\vskip 0.8 cm
\centerline{\BIGtitle Extension theorems for reductive group schemes}

\vskip 0.8 cm
\centerline{\bigtitle Adrian Vasiu}

\vskip 0.4 cm
\centerline{December 11, 2015, accepted in final form for publication in J. Algebra \& Number Theory}
\footline={\hfill}
\vskip 0.4 cm
\noindent
{\bf ABSTRACT.} We prove several basic extension theorems for reductive group schemes via extending Lie algebras and via taking schematic closures. We also prove that for each scheme $Y$, the category in groupoids of adjoint group scheme over $Y$ whose Lie algebras $\scrO_Y$-modules have perfect Killing forms is isomorphic via the differential functor to the category in groupoids of Lie algebras $\scrO_Y$-modules which have perfect Killing forms and which as $\scrO_Y$-modules are coherent and locally free.
  

\footline={\hss\tenrm \folio\hss}
\pageno=1

\bigskip
\noindent
{\boldsectionfont 1. Introduction}
\bigskip

A group scheme $H$ over a scheme $S$ is called {\it reductive} if the morphism $H\to S$ has the following two properties: (i) it is smooth and affine (and therefore of finite presentation) and (ii) its geometric fibres are reductive groups over spectra of fields and therefore are connected (cf. [DG, Vol. III, Exp. XIX, Subsects. 2.7, 2.1, and 2.9]). If moreover the {\it center} of $H$ is trivial, then $H$ is called an {\it adjoint} group scheme over $S$. Let $\scrO_S$ be the structure ring sheaf of $S$. Let $Lie(H)$ be the {\it Lie algebra} $\scrO_S$-module of $H$. As a $\scrO_S$-module, $Lie(H)$ is coherent and locally free. 

The main goal of the paper is to prove Theorems 1.2 and 1.4 below (see Sections 3 and 4) and to apply them and their proofs for obtaining new extension theorems for homomorphisms between reductive group schemes (see Section 5). We begin by introducing two groupoids on sets (i.e., two categories whose morphisms are all isomorphisms).

\medskip\smallskip\noindent
{\bf 1.1. Two groupoids on sets.} Let $Y$ be an arbitrary scheme. Let $\text{Adj-perf}_Y$ be the category whose objects are adjoint group schemes over $Y$ with the property that their Lie algebras $\scrO_Y$-modules have perfect {\it Killing forms} (i.e., the Killing forms induce naturally $\scrO_Y$-linear isomorphisms from them into their duals) and whose morphisms are isomorphisms of group schemes. Let $\text{Lie-perf}_Y$ be the category whose objects are Lie algebras $\scrO_Y$-modules which have perfect Killing forms and which as $\scrO_Y$-modules are coherent and locally free and whose morphisms are isomorphisms of Lie algebras $\scrO_Y$-modules.

\medskip\smallskip\noindent
{\bf 1.2. Theorem.} {\it Let $\scrL_Y:\text{Adj-perf}_Y\to\text{Lie-perf}_Y$ be the functor which associates to a morphism $f:G\arrowsim H$ of $\text{Adj-perf}_Y$ the morphism $df:Lie(G)\arrowsim Lie(H)$ of $\text{Lie-perf}_Y$ which is the differential of $f$. Then the functor $\scrL_Y$ is an equivalence of categories.}

\medskip
We have a variant of this theorem for simply connected semisimple group schemes instead of adjoint group schemes, cf. Corollary 3.8. This theorem implies the classification of Lie algebras over fields of characteristic at least $3$ that have non-degenerate Killing forms obtained previously by Curtis, Seligman, Mills, Block--Zassenhaus, and Brown in [C], [MS], [Mi], [BR], [S], and [Br] (see Remark 3.7 (a)). The functor $\scrL_Y$ is an equivalence of non-empty categories if and only if $Y$ is a $\Spec Z[{1\over 2}]$-scheme, cf. Corollary 3.9. Directly from the Theorem 1.2 we get our first extension result:

\medskip\smallskip\noindent
{\bf 1.3. Corollary.} {\it We assume that $Y=\Spec A$ is an affine scheme. Let $K$ be the ring of fractions of $A$. Let $G_K$ be an adjoint group scheme over $\Spec K$ such that the symmetric bilinear Killing form on the Lie algebra $\Lie(G_K)$ of $G_K$ is perfect (i.e., it induces naturally a $K$-linear isomorphism $\Lie(G_K)\arrowsim\Hom_K(\Lie(G_K),K)$). We assume that there exists a Lie algebra $\grg$ over $A$ such that the following two properties hold:

\medskip
{\bf (i)} we have an identity $\Lie(G_K)=\grg\otimes_A K$ and the $A$-module $\grg$ is projective and finitely generated;

\smallskip
{\bf (ii)} the symmetric bilinear Killing form on $\grg$ is perfect. 

\medskip
Then there exists a unique adjoint group scheme $G$ over $Y$ which extends $G_K$ and such that we have an identity $\Lie(G)=\grg$ that extends the identity of the property (i).}

\medskip
Let $U$ be an open, Zariski dense subscheme of $Y$. We call the pair $(Y,Y\setminus U)$ {\it quasi-pure} if each finite \'etale cover of $U$ extends uniquely to a finite \'etale cover of $Y$ (to be compared with [G2, Exp. X, Def. 3.1]).

\medskip\smallskip\noindent
{\bf 1.4. Theorem.} {\it We assume that $Y$ is a normal, noetherian scheme and the codimension of $Y\setminus U$ in $Y$ is at least $2$. Then the following two properties hold:

\medskip
{\bf (a)} Let $G_U$ be an adjoint group scheme over $U$. We assume that the Lie algebra $\scrO_U$-module $Lie(G_U)$ of $G_U$ extends to a Lie algebra $\scrO_Y$-module that is a locally free $\scrO_Y$-module. Then $G_U$ extends uniquely to an adjoint group scheme $G$ over $Y$.

\medskip
{\bf (b)} Let $H_U$ be a reductive group scheme over $U$. We assume that the pair $(Y,Y\setminus U)$ is quasi-pure and that the Lie algebra $\scrO_U$-module $Lie(G_U)$ of the adjoint group scheme $G_U$ of $H_U$ extends to a Lie algebra $\scrO_Y$-module that is a locally free $\scrO_Y$-module. Then $H_U$ extends uniquely to a reductive group scheme $H$ over $Y$.}

\medskip
The proof of the Theorem 1.2 we include combines the cohomology theory of Lie algebras with a simplified variant of [V1, Claim 2, p. 464] (see Subsections 3.3 and 3.5). The proof of Theorem 1.4 (a) is an application of [CTS, Cor. 6.12] (see Subsection 4.1). The classical purity theorem of Nagata and Zariski (see [G2, Exp. X, Thm. 3.4 (i)]) says that the pair $(Y,Y\setminus U)$ is quasi-pure, provided $Y$ is {\it regular} and $U$ contains all points of $Y$ of codimension $1$ in $Y$. In such a case, a slightly weaker form of Theorem 1.4 (b) was obtained in [CTS, Thm. 6.13]. In general, the hypotheses of the Theorem 1.4 are needed (see Remarks 4.3). See [MB] (resp. [FC], [V1], [V2], and [VZ]) for different analogues of Theorems 1.4 (a) and (b) for Jacobian (resp. for abelian) schemes. For instance, in [VZ, Cor. 1.5] it is proved that if $Y$ is a regular, formally smooth scheme over the spectrum of a discrete valuation ring of mixed characteristic $(0,p)$ and index of ramification at most $p-1$ and if $U$ contains all points of $Y$ that are of either characteristic $0$ or codimension $1$ in $Y$, then each abelian scheme over $U$ extends uniquely to an abelian scheme over $Y$. 

Section 2 presents notations and basic results. In Section 3 we prove Theorem 1.2. In Section 4 we prove Theorem 1.4.

Section 5 contains two results on extending homomorphisms between reductive group schemes. The first one is an application of Theorem 1.4 (b) (see Proposition 5.1) and pertains to extensions of homomorphisms in codimension at least 2 over normal bases. The second one pertains to extensions of homomorphisms via schematic closures and refines [V1, Lemma 3.1.6] (see Proposition 5.2); its role is to achieve natural reductions such as the reduction to the case of either a torus or a semisimple group scheme. 

Our main motivation for Theorems 1.2 and 1.4 stems from the meaningful applications to crystalline cohomology one gets by combining them with either Faltings' results of [F, Sect. 4] (see [V1] and [V7]) or de Jong's extension theorem [dJ, Thm. 1.1] (see  [V6]). The manuscript [V6] applies this paper to extend our prior work on integral canonical models of Shimura varieties of Hodge type in unramified mixed characteristic $(0,p)$ with $p\ge 5$ (see [V1]), to unramified mixed characteristics $(0,2)$ and $(0,3)$. In addition, this paper can be used to get relevant simplifications of certain parts of the mentioned prior work (see [V7]).

\bigskip
\noindent
{\boldsectionfont 2. Preliminaries}
\bigskip 

Our notations are gathered in Subsection 2.1. In Subsections 2.2 to 2.5 we include four basic results that are of different nature and that are often used in Sections 3 to 5.

\medskip\smallskip\noindent
{\bf 2.1. Notations and conventions.} 
If $K$ is a field, let $\bar K$ be an algebraic closure of it. Let $H$ be a reductive group scheme over a scheme $S$. Let $Z(H)$, $H^{\der}$, $H^{\ad}$, and  $H^{\ab}$ denote the center, the {\it derived group} scheme, the adjoint group scheme, and the {\it abelianization} of $H$ (respectively). We have $H^{\ab}=H/H^{\der}$ and $H^{\ad}=H/Z(H)$. The center $Z(H)$ is a group scheme of {\it multiplicative type}, cf. [DG, Vol. III, Exp. XXII, Cor. 4.1.7]. Let $Z^0(H)$ be the maximal {\it torus} of $Z(H)$; the {\it quotient} group scheme $Z(H)/Z^0(H)$ is a finite, flat group scheme over $S$ of multiplicative type. Let $H^{\sc}$ be the simply connected semisimple group scheme cover of the derived group scheme $H^{\der}$. 

See [DG, Vol. III, Exp. XXII, Cor. 4.3.2] for the quotient group scheme $H/F$ of $H$ by a flat, closed subgroup scheme $F$ of $Z(H)$ which is of multiplicative type. 

If $X$ or $X_S$ is an $S$-scheme, let $X_{A_1}$ (resp. $X_{S_1}$) be its pull back via a morphism $\Spec A_1\to S$ (resp. $S_1\to S$). 

If $S$ is either affine or integral, let $K_S$ be the ring of fractions of $S$. If $S$ is a normal, noetherian, integral scheme, let $\scrD(S)$ be the set of local rings of $S$ that are discrete valuation rings. 

Let $\dbG_{m,S}$ be the rank 1 split torus over $S$; similarly, the group schemes $\dbG_{a,S}$, $\pmb{GL}_{d,S}$ with $d\in\db N^*$, etc., will be understood to be over $S$. Let $Lie(H)$ be the Lie algebra $\scrO_S$-module of $H$. If $S=\Spec A$ is affine, then let $\dbG_{m,A}:=\dbG_{m,S}$, etc., and let $\Lie(F)$ be the Lie algebra over $A$ of a closed subgroup scheme $F$ of $H$. As $A$-modules, we identify $\Lie(F)=\Ker(F(A[x]/x^2)\to F(A))$, where the $A$-epimorphism $A[x]/(x^2)\twoheadrightarrow A$ takes $x$ to $0$. The Lie bracket on $\Lie(F)$ is defined by taking the (total) differential of the commutator morphism $[\, ,]:F\times_S F\to F$ at identity sections. If $S=\Spec A$ is affine, then $\Lie(H)=Lie(H)(S)$ is the Lie algebra over $A$ of global sections of $Lie(H)$ and it is a projective, finitely generated $A$-module. 

If $N$ is a projective, finitely generated $A$-module, let $N^*:=\Hom_A(N,A)$, let $\pmb{GL}_N$ be the reductive group scheme over $\Spec A$ of linear automorphisms of $N$, and let $\gl_N:=\Lie(\pmb{GL}_N)$. Thus $\gl_N$ is the Lie algebra associated to the $A$-algebra $\End_A(N)$. 
A bilinear form $b_N:N\times N\to A$ on $N$ is called {\it perfect} if it induces an $A$-linear map $N\to N^*$ that is an isomorphism. If $b_N$ is symmetric, then its {\it kernel} is the $A$-submodule 
$$\Ker(b_N):=\{a\in N|b_N(a,b)=0\;\forall b\in N\}$$ 
of $N$. For a Lie algebra $\grg$ 
over $A$ that is a projective, finitely generated $A$-module, let $\text{ad}:\grg\to\gl_{\grg}$  be the adjoint representation of $\grg$ and let $\scrK_{\grg}:\grg\times\grg\to A$ be the Killing form on $\grg$. For $a$, $b\in\grg$ we have $\text{ad}(a)(b)=[a,b]$ and $\scrK_{\grg}(a,b)$ is the trace of the endomorphism $\ad(a)\circ\ad(b)$ of $\grg$. The kernel $\Ker(\scrK_{\grg})$ is an ideal of $\grg$.

We denote by $k$ an arbitrary field. Let $n\in\dbN^*$. See [B2, Ch. VI, Sect. 4] and [H1, Ch. III, Sect. 11] for the classification of connected Dynkin diagrams. For 
$$\flat\in\{A_n,B_n,C_n|n\in\dbN^*\}\cup\{D_n|n\Ge 3\}\cup\{E_6,E_7,E_8,F_4,G_2\}$$ 
we say that $H$ is of {\it isotypic $\flat$ Dynkin type} if the connected Dynkin diagram of each simple factor of an arbitrary geometric fibre of $H^{\ad}$, is $\flat$; if $H^{\ad}$ is absolutely simple we drop the word isotypic. We recall that $A_1=B_1=C_1$, $B_2=C_2$, and $A_3=D_3$.

\medskip\smallskip\noindent
{\bf 2.2. Proposition.} {\it Let $Y$ be a normal, noetherian, integral scheme. Let $K:=K_Y$. 

\medskip
{\bf (a)} If $Y=\Spec A$ is affine, then inside the field $K$ we have $A=\cap_{V\in\scrD(Y)} V$.

\smallskip
{\bf (b)} Let $U$ be an open subscheme of $Y$ such that $Y\setminus U$ has codimension in $Y$ at least $2$. Let $W$ be an affine $Y$-scheme of finite type. Then the natural restriction map $\Hom_Y(Y,W)\to \Hom_Y(U,W)$ is a bijection. If moreover $W$ is integral, normal and such that we have $\scrD(W)=\scrD(W_U)$, then $W$ is determined (up to unique isomorphism) by $W_U$. 

\smallskip
{\bf (c)} Suppose that $Y=\Spec A$ is local, regular, and has dimension $2$. Let $y$ be the closed point of $Y$ and let $U:=Y\setminus\{y\}$. Then each locally free $\scrO_U$-module of finite rank, extends uniquely to a free $\scrO_Y$-module.}

\medskip
\proof
See [M, (17.H), Thm. 38] for (a). To check (b), we can assume $Y=\Spec A$ is affine. We write $W=\Spec B$. The $A$-algebra of global functions of $U$ is $A$, cf. (a). We have $\Hom_Y(U,W)=\Hom_A(B,A)=\Hom_Y(Y,W)$. If moreover $B$ is a normal ring and we have $\scrD(W)=\scrD(W_U)$, then $B$ is uniquely determined by $\scrD(W_U)$ (cf. (a)) and therefore by $W_U$. From this (b) follows. See [G2, Exp. X, Lemma 3.5] for (c).\endproof

\medskip\smallskip\noindent
{\bf 2.3. Proposition.} {\it Let $G$ be a reductive group scheme over a scheme $Y$. Then the functor on the category of $Y$-schemes that parametrizes maximal tori of $G$, is representable by a smooth, separated $Y$-scheme of finite type. Thus locally in the \'etale topology of $Y$, $G$ has split, maximal tori.}

\medskip
\proof
See [DG, Vol. II, Exp. XII, Cor. 1.10] for the first part. 
The second part follows easily from the first part (see also [DG, Vol. III, Exp. XIX, Prop. 6.1]).\endproof

\medskip\noindent
{\bf 2.3.1.  Lemma.} {\it Let $Y$ be a reduced scheme. Let $G$ be a reductive group scheme over $Y$. Let $K:=K_Y$. Let $f_K:G_K'\to G_K$ be a central isogeny of reductive group schemes over $\Spec K$. We assume that either $G$ is split or $Y$ is normal. We have:

\medskip
{\bf (a)} There exists (up to a canonical identification) at most one central isogeny $f:G'\to G$ that extends $f_K:G_K'\to G_K$. If $Y$ is integral (i.e., $K$ is a field), then there exists a unique central isogeny $f:G'\to G$ that extends $f_K:G_K'\to G_K$. 

\smallskip
{\bf (b)} If $Y$ is normal and integral, then $G'$ is the normalization of $G$ in (the field of fractions of) $G_K'$.}

\medskip
\proof
We first prove (a) in the case when $G$ is split. Let $T$ be a split, maximal torus of $G$. We first prove the existence part: thus $K$ is a field. As $f_K$ is a central isogeny, the inverse image $T_K'$ of $T_K$ in $G_K'$ is a split torus. Thus $G_K'$ is split. Let $\scrR'\to\scrR$ be the $1$-morphism of root data in the sense of [DG, Vol. III, Exp. XXI, Def. 6.8.1] which is associated to the central isogeny $f_K:G_K'\to G_K$ that extends the isogeny $T_K'\to T_K$. Let $\tilde f:\tilde G'\to G$ be a central isogeny of split, reductive group schemes over $Y$ which extends an isogeny of split tori $\tilde T'\to T$ and for which the $1$-morphism of root data associated to it and to the isogeny $\tilde T'\to T$, is $\scrR'\to\scrR$ (cf. [DG, Vol. III, Exp. XXV, Thm. 1.1]). From loc. cit. we also get that there exists an isomorphism $i_K:\tilde G_K'\arrowsim G_K'$ such that we have $\tilde f_K=f_K\circ i_K$. Obviously, $i_K$ is unique. Let $G'$ be the unique group scheme over $Y$ such that $i_K$ extends (uniquely) to an isomorphism $i:\tilde G'\arrowsim G'$. Let $f:=\tilde f\circ i^{-1}:G'\to G$; it is a central isogeny. 

To check the uniqueness part, we consider two central isogenies $G'\to G$ and $G'_1\to G$ that extend a central isogeny $f_K:G'_K\to G_K$ (thus $G'_K=G'_{1,K}$). Let $G'_2$ be the schematic closure of $G_K'$ embedded diagonally into the product $G'\times_Y G'_1$. We are left to check that the two projections $\pi_1:G'_2\to G'$ and $\pi_2:G'_2\to G'_1$ are isomorphisms as in such a case the composite isomorphism $\pi_2\circ\pi_1^{-1}:G'\arrowsim G_1'$ is an isomorphism that extends the identity $G'_K=G'_{1,K}$. This statement is local for the \'etale topology of $Y$ and therefore we can assume based on Proposition 2.3 that the inverse images of $T$ to $G'$ and $G_1'$ are split tori. From this and [DG, Vol. III, Exp. XXIII, Thm. 4.1] we get that there exists a unique isomorphism $\theta:G_1\arrowsim G'_1$ which extends the identity $G'_K=G_K$. This implies that $G_2'$ is the graph of $\theta$ and therefore the two projections $\pi_1$ and $\pi_2$ are isomorphisms.
We conclude that (a) holds if $G$ is split.

We now prove simultaneously (a) and (b) in the case when $Y$ is normal. If a $G'$ as in (a) exists, then it is a smooth scheme over the normal scheme $Y$ and thus it is a normal scheme; from this and the fact that $f:G'\to G$ is a finite morphism, we get that $G'$ is the normalization of $G$ in $G_K'$ and in particular it is unique. 

Thus to end the proof of the Lemma, it suffices to show that the normalization $G'$ of $G$ in $G_K'$ is a reductive group scheme equipped with a central isogeny $f:G'\to G$. This is a local statement for the \'etale topology of $Y$. As each connected, \'etale  scheme over $Y$ is a normal, integral scheme, based on Proposition 2.3 we can assume that $G$ has a split, maximal torus $T$. Thus the fact that $G'$ is a reductive group scheme equipped with a central isogeny $f:G'\to G$ follows from the previous three paragraphs.\endproof

\medskip\noindent
{\bf 2.3.2. Lemma.} {\it Let $Y=\Spec A$ be an affine scheme. Let $K:=K_Y$. Let $T$ be a torus over $Y$ equipped with a homomorphism $\rho:T\to G$, where $G$ is a reductive group scheme over $Y$. Then the following three properties hold:

\medskip
{\bf (a)} the kernel $\Ker(\rho)$ is a group scheme over $Y$ of multiplicative type; 

\smallskip
{\bf (b)}  the kernel $\Ker(\rho)$ is trivial (resp. finite) if and only if the kernel $\Ker(\rho_K)$ is  trivial (resp. is finite);

\smallskip
{\bf (c)} the quotient group scheme $T/\Ker(\rho)$ is a torus and we have a closed embedding homomorphism $T/\Ker(\rho)\hookrightarrow G$.}

\medskip
\proof
The statements of the Lemma are local for the \'etale topology of $Y$. Thus we can assume that $Y$ is local and (cf. Proposition 2.3) that $T$ and $G$ are split. As $Y$ is connected, the split reductive group scheme $G$ has constant Lie type. Thus $G$ is the pull back to $Y$ of a reductive group scheme $G_{\dbZ}$ over $\Spec\dbZ$, cf. [DG, Vol. III, Exp. XXV, Cor. 1.2]. As $G_{\dbZ}$ can be embedded into a general linear group scheme over $\Spec\dbZ$ (for instance, cf. [DG, Vol. I, Exp. VI${}_B$, Rm. 11.11.1]), there exists a closed embedding homomorphism $G\hookrightarrow \pmb{GL}_M$, with $M$ a free $A$-module of rank $d\in\db N^*$. By replacing $\rho$ with its composite with this closed embedding homomorphism $G\hookrightarrow \pmb{GL}_M$, we can assume that $G=\pmb{GL}_M$ is a general linear group scheme over $Y$. The representation of $T$ on $M$ is a finite direct sum of representations of $T$ of rank $1$, cf. [J, Part I, Subsect. 2.11]. Thus $\rho$ factors as the composite of a homomorphism $\rho_1:T\to\dbG_{m,A}^m$ with a closed embedding homomorphism $\dbG_{d,A}^m\hookrightarrow \pmb{GL}_M$. The kernel $\Ker(\rho_1)$ is a group scheme over $Y$ of multiplicative type, cf. [DG, Vol. II, Exp. IX, Prop. 2.7 (i)]. As $\Ker(\rho)=\Ker(\rho_1)$, we get that (a) holds. As (a) holds,  $\Ker(\rho)$ is flat over $Y$ as well as the extension of a finite, flat group scheme $T_1$ by a torus $T_0$. But $T_1$ (resp. $T_0$) is a trivial group scheme if and only if $T_{1,K}$ (resp. $T_{0,K}$) is trivial. From this (b) follows. The quotient group scheme $T/\Ker(\rho)$ exists and is a closed subgroup scheme of $\dbG_{m,A}^m$ that is of multiplicative type, cf. [DG, Vol. II, Exp. IX, Prop. 2.7 (i) and Cor. 2.5]. As the fibres of $T/\Ker(\rho)$ are tori, we get that $T/\Ker(\rho)$ is a torus. Thus (c) holds.\endproof

\medskip\smallskip\noindent
{\bf 2.4. Lemma.} {\it Suppose that $k=\bar k$. Let $H$ be a reductive group over $\Spec k$. Let $\grn$ be a non-zero ideal of $\Lie(H)$ which is a simple left $H$-module. We assume that there exists a maximal torus $T$ of $H$ such that we have $\Lie(T)\cap\grn=0$. Then $\char(k)=2$ and $H^{\der}$ has a normal, subgroup $F$ which is isomorphic to $\pmb{SO}_{2n+1,k}$ for some $n\in\dbN^*$ and for which we have an inclusion $\grn\subseteq\Lie(F)$.} 

\medskip
\proof
This is only a variant of [V4, Lemma 2.1].\endproof

\medskip\noindent
{\bf 2.4.1. Remark.} If $\grn $ is assumed to be a restricted Lie subalgebra of $\Lie(H)$ (for instance, this holds if $\grn$ is the Lie algebra of a subgroup of $H$), then there exists a purely inseparable isogeny $H\to H/\grn$ (see [Bo, Ch. V, Prop. 7.4]) and in this case Lemma 2.4 can be also deduced easily from [PY, Lemma 2.2] applied to such isogenies with $H^{\ad}$ absolutely simple. In this paper, Lemma 2.4 will be applied only in such situations in which $\grn$ is the Lie algebra of a subgroup of $H$.

\medskip\smallskip\noindent
{\bf 2.5. Theorem.} {\it Let $f:G_1\to G_2$ be a homomorphism between group schemes over a scheme $Y$. We assume that $G_1$ is reductive, that $G_2$ is separated and of finite presentation, and that all fibres of $f$ are closed embeddings. Then $f$ is a closed embedding.} 

\medskip
\proof
As $G_1$ is of finite presentation over $Y$, the homomorphism $f$ is locally of finite type. As the fibres of $f$ are closed embeddings and thus monomorphisms, $f$ itself is a monomorphism (cf. [DG, Vol. I, Exp. VI${}_B$, Cor. 2.11]). Thus the Theorem follows from [DG, Vol. II, Exp. XVI, Cor. 1.5 a)].\endproof  

\medskip\noindent
{\bf 2.5.1. Lemma.} {\it Let $G$ be an adjoint group scheme over an affine scheme $Y=\Spec A$. Let $\Aut(G)$ be the group scheme over $Y$ of automorphisms of $G$. Then the natural adjoint representation $\Ad:\Aut(G)\to \pmb{GL}_{\Lie(G)}$ is a closed embedding.}

\medskip
\proof
To prove the Lemma, we can work locally in the \'etale topology of $Y$ and therefore (cf. Proposition 2.3) we can assume that $G$ is split and that $Y$ is connected. We have a short exact sequence $1\to G\to\Aut(G)\to C\to 1$ that splits (cf. [DG, Vol. III, Exp. XXIV, Thm. 1.3]), where $C$ is a finite, \'etale, constant group scheme over $Y$. Thus $G$ is the identity component of $\Aut(G)$ and $\Aut(G)$ is a finite disjoint union of right translates of $G$ via certain $Y$-valued points of $\Aut(G)$. If the fibres of $\Ad$ are closed embeddings, then the restriction of $\Ad$ to $G$ is a closed embedding (cf. Theorem 2.5) and thus also the restriction of $\Ad$ to any right translate of $G$ via a $Y$-valued point of $\Aut(G)$ is a closed embedding. The last two sentences imply that $\Ad$ is a closed embedding. Thus to end the proof, we are left to check that the fibres of $\Ad$ are closed embeddings. For this, we can assume that $A$ is an algebraically closed field. 

As $G$ is adjoint and $A$ is a field, the restriction of $\Ad$ to $G$ is a closed embedding. Thus the representation $\Ad$ is a closed embedding if and only if each element $g\in\Aut(G)(A)$ that acts trivially on $\Lie(G)$, is trivial. We show that the assumption that there exists a non-trivial such element $g$ leads to a contradiction. For this, we can assume that $G$ is absolutely simple and that $g$ is a non-trivial outer automorphism of $G$. Let $T$ be a maximal torus of a Borel subgroup $B$ of $G$ and let $n$ be the dimension of $T$. 

For $t\in\Lie(T)$, let $C_G(t)$ be its centralizer in $G$. It is a subgroup of $G$ that contains $T$. In this paragraph we check that, as $G$ is adjoint, we can choose $t$ such that we have $C_G(t)^0=T$. We consider the root decomposition $\Lie(G)=\Lie(T)\bigoplus_{\alpha\in\Phi} \grg_{\alpha}$ with respect to $T$, where $\Phi$ is the root system of $G$ and where each $\grg_{\alpha}$ is a one dimensional $A$-vector space normalized by $T$. Let $\Delta$ be the basis for $\Phi$ such that we have $\Lie(B)=\Lie(T)\bigoplus_{\alpha\in\Delta} \grg_{\alpha}$. As $G$ is adjoint, $\Delta$ is a basis for the dual $A$-vector space $\Lie(T)^*$ (to be compared with [DG, Vol. III, Exp. XXI, Def. 6.2.6 and Exp. XXII, Def. 4.3.3]). Thus for each root $\alpha\in\Delta$, $\Ker(\alpha)$ is a $A$-vector subspace of $\Lie(T)$ of dimension $n-1$. As each $\alpha\in\Phi$ is conjugate under the Weyl group of $\Phi$ (equivalently of $G$) to an element of $\Delta$ (see [H1, Ch. III, Sect. 10, Thm.]), we get that for each $\alpha\in\Delta$ its kernel $\Ker(\alpha)$ is a $A$-vector subspace of $\Lie(T)$ of dimension $n-1$. We choose $t\in\Lie(T)\setminus (\cup_{\alpha\in\Phi} \Ker(\alpha))$. This implies that we have $\Lie(C_G(t))=\Lie(T)$. From this and the fact that $T$ is a subgroup of $C_G(t)$, we get that $C_G(t)$ is a smooth group of dimension $n$ and therefore that $C_G(t)^0=T$. 

As $g$ fixes $t$ and $\Lie(B)$, $g$ normalizes both $C_G(t)^0=T$ and $B$. But it is well known that a non-trivial outer automorphism $g$ of $G$ that normalizes both $T$ and $B$, can not fix $\Lie(B)$. Contradiction. Thus $\Ad$ is a closed embedding.\endproof

\medskip
We follow the ideas of [V1, Prop. 3.1.2.1 c) and Rm. 3.1.2.2 3)] in order to prove the next Proposition.

\medskip\noindent
{\bf 2.5.2. Proposition.} {\it Let $V$ be a discrete valuation ring whose residue field is $k$. Let $Y=\Spec V$ and let $K:=K_Y$. Let $f:H_1\to H_2$ be a homomorphism between flat, finite type, affine group schemes over $Y$ such that $H_1$ is a reductive group scheme and the generic fibre $f_K:H_{1,K}\to H_{2,K}$ of $f$ is a closed embedding. We have:

\medskip
{\bf (a)} The subgroup scheme $\Ker(f_k:H_{1,k}\to H_{2,k})$ of $H_{1,k}$ has a trivial intersection with each torus $T_{1,k}$ of $H_{1,k}$. In particular, we have $\Lie(\Ker(f_k))\cap\Lie(T_{1,k})=0$.

\smallskip
{\bf (b)} The homomorphism $f$ is finite. 

\smallskip
{\bf (c)} If $\char(k)=2$, we assume that $H_{1,K}$ has no normal subgroup that is adjoint of isotypic $B_n$ Dynkin type for some $n\in\dbN^*$. Then $f$ is a closed embedding.}

\medskip
\proof
Let $\rho:H_2\hookrightarrow \pmb{GL}_M$ be a closed embedding homomorphism, with $M$ a free $V$-module of finite rank (cf. [DG, Vol. I, Exp. VI${}_B$, Rm. 11.11.1]). To prove the Proposition we can assume that $V$ is complete, that $k=\bar k$, and that $f_K:H_{1,K}\to H_{2,K}$ is an isomorphism. Let $H_{0,k}:=\Ker(f_k)$. We now show that the group scheme $H_{0,k}\cap T_{1,k}$ is trivial by adapting arguments from [V1, Rm. 3.1.2.2 3) and proof of Lemma 3.1.6]. As $V$ is strictly henselian, the maximal torus $T_{1,k}$ of $H_{1,k}$ is split and (cf. Proposition 2.3) it lifts to a maximal torus $T_1$ of $H_1$. The restriction of $\rho\circ f$ to $T_1$ has a trivial kernel (as its fibre over $\Spec K$ is trivial, cf. Lemma 2.3.2 (b)) and therefore it is a closed embedding (cf. Lemma 2.3.2 (c)). Thus the restriction of $f$ to $T_1$ is a closed embedding homomorphism $T_1\hookrightarrow H_2$. Therefore the intersection $H_{0,k}\cap T_{1,k}$ is a trivial group scheme. Thus (a) holds.

We check (b). The identity component of the reduced scheme of $\Ker(f_k)$ is a reductive group that has $0$ rank (cf. (a)) and therefore it is a trivial group. Thus $f$ is a quasi-finite, birational morphism. From Zariski Main Theorem (see [G1, Thm. (8.12.6)]) we get that $H_1$ is an open subscheme of the normalization $H_2^{\n}$ of $H_2$. Let $H_3$ be the smooth locus of $H_2^{\n}$ over $\Spec V$; it is an open subscheme of $H_2^{\n}$ that contains $H_1$. As $H_3$ is an open subscheme of the affine scheme $H_2^{\n}$, it is a quasi-affine scheme. 

As $H_3$ is smooth over $\Spec V$, the products $H_3\times_{\Spec V} H_2^{\n}$ and $H_2^{\n}\times_{\Spec V} H_3$ are smooth over $H_2^{\n}$ and thus are normal schemes. The product $H_2^{\n}\times_{\Spec V} H_2^{\n}$ is a flat scheme over $\Spec V$ whose generic fibre is smooth over $\Spec K$. Its normalization $(H_2^{\n}\times_{\Spec V} H_2^{\n})^{\n}$ contains both $H_3\times_{\Spec V} H_2^{\n}$ and $H_2^{\n}\times_{\Spec V} H_3$ as open subschemes and is equipped with a finite surjective morphism $(H_2^{\n}\times_{\Spec V} H_2^{\n})^{\n}\to H_2^{\n}\times_{\Spec V} H_2^{\n}$ whose generic fibre is an isomorphism. The product morphism $H_2\times_{\Spec V} H_2\to H_2$ induces a natural product type of morphism $\Theta:(H_2^{\n}\times_{\Spec V} H_2^{\n})^{\n}\to H_2^{\n}$. Its restrictions to $H_3\times_{\Spec V} H_2^{\n}$ and $H_2^{\n}\times_{\Spec V} H_3$ induce product type of morphisms $H_3\times_{\Spec V} H_2^{\n}\to H_2^{\n}$ and $H_2^{\n}\times_{\Spec V} H_3\to H_2^{\n}$. This implies that for each valued point $z\in H_2^{\n}(V)$ it makes sense to speak about the left $zH_3$ and the right $H_3z$ translations of $H_3$ through $z$; they are smooth open subschemes of $H_2^{\n}$ and thus of $H_3$. This implies that $H_3(V)=H_2^{\n}(V)$ and that $\Theta$ restricts to a product morphism $H_3\times_{\Spec V} H_3\to H_3$. The inverse automorphisms of the $\Spec V$-schemes $H_1$ and $H_2$ induce an inverse automorphism of the $\Spec V$-scheme $H_2^{\n}$ which restricts to an inverse automorphism of the $\Spec V$-scheme $H_3$. With respect to its product morphism, its inverse automorphism, and its identity section inherited from $H_1$, $H_3$ gets the structure of a (quasi-affine) group scheme over $\Spec V$ that is of finite type. 

As $V$ is complete, it is also excellent (cf. [M, Sect. 34]). Thus the morphism $H_2^{\n}\to H_2$ is finite. The homomorphism $f$ is finite if and only if $H_1=H_2^{\n}$ and thus if and only if the set $H_2^{\n}(k)\setminus H_1(k)$ is empty. We show that the assumption that $H_1\neq H_2^{\n}$ leads to a contradiction. Let $x\in H_2^{\n}(k)\setminus H_1(k)$. From [V5, Lemma 4.1.5] applied to the completion of the local ring of $x$ in $H_2^{\n}$, we get that there exists a finite, flat discrete valuation ring extension $V'$ of $V$ for which we have a valued point $z'\in H^{\n}_2(V')$ that lifts $x$ (we recall that loc. cit. is only a local version of the global result [G1, Cor. (17.16.2)]). The flat $\Spec V'$-scheme $H^{\n}_{2,V'}$ might not be normal but we have $H_1\neq H_2^{\n}$ if and only if $H_{1,V'}\neq H_{2,V'}^{\n}$. Thus to reach a contradiction we can replace $V$ by $V'$ and therefore we can assume that $x$ is such that there exists a valued point $z\in H^{\n}_2(V)=H_3(V)$ which lifts $x$. As $x\in H_2^{\n}(k)\setminus H_1(k)$, we have $z\in H_3(V)\setminus H_1(V)$. As $H_1$ is a subgroup scheme of $H_3$, all fibers of the homomorphism $H_1\to H_3$ are closed. From this and Theorem 2.5 we get that $H_1$ is a closed subscheme of $H_3$. Thus, as $H_3$ is an integral scheme and as $H_{3,K}=H_{1,K}$, we get that $H_1=H_3$. This contradicts the fact that $z\in H_3(V)\setminus H_1(V)$. Thus (b) holds. 

We check (c). We show that the assumption that $\Lie(H_{0,k})\neq 0$ leads to a contradiction. From Lemma 2.4 applied to $H_{1,k}$ and to any simple $H_{1,k}$-submodule of the left $H_{1,k}$-module $\Lie(H_{0,k})$, we get that $\char(k)=2$ and that $H_{1,k}$ has a normal subgroup $H_{4,k}$ isomorphic to $\pmb{SO}_{2n+1,k}$ for some $n\in\dbN^*$. As $H_{4,k}$ is adjoint, we have a product decomposition $H_{1,k}=H_{4,k}\times_{\Spec k} H_{5,k}$ of reductive groups. It lifts (cf. [DG, Vol. III, Exp. XXIV, Prop. 1.21]) to a product decomposition $H_1=H_4\times_{\Spec V} H_5$, where $H_4$ is isomorphic to $\pmb{SO}_{2n+1,V}$ and where $H_5$ is a reductive group scheme over $\Spec V$. This contradicts the extra hypothesis of (c). Thus we have $\Lie(H_{0,k})=0$. Therefore  $H_{0,k}$ is a finite, \'etale, normal subgroup of $H_{1,k}$. But $H_{1,k}$ is connected and thus its action on $H_{0,k}$ via inner conjugation is trivial. Therefore we have $H_{0,k}\leqslant Z(H_1)_k\leqslant T_{1,k}$. Thus $H_{0,k}=H_{0,k}\cap T_{1,k}$ is the trivial group, cf. (a). In other words, the homomorphism $f_k:H_{1,k}\to H_{2,k}$ is a closed embedding. Thus $f:H_1\to H_2$ is a closed embedding homomorphism, cf. Theorem 2.5.\endproof

\medskip\noindent
{\bf 2.5.3. Remark.} See [V3, Thm. 1.2 (b)] and [PY, Thm. 1.2] for two other proofs of Proposition 2.5.2 (c).

\bigskip
\noindent
{\boldsectionfont 3. Lie algebras with perfect Killing forms}
\bigskip 

Let $A$ be a commutative $\dbZ$-algebra. Let $\grg$ be a Lie algebra over $A$ which as an $A$-module is projective and finitely generated. In this Section we will assume that the Killing form $\scrK_{\grg}$ on $\grg$ is perfect. Let $U_{\grg}$ be the enveloping algebra of $\grg$ i.e., the quotient of the tensor algebra $T_{\grg}$ of $\grg$ by the two-sided ideal of  $T_{\grg}$ generated by the subset $\{x\otimes y-y\otimes x-[x,y]|x,y\in\grg\}$ of $T_{\grg}$. Let $Z(U_{\grg})$ be the center of $U_{\grg}$. The categories of left $\grg$-modules and of left $U_{\grg}$-modules are canonically identified. We view $\grg$ as a left $\grg$-module via the adjoint representation $\ad:\grg\to\grg\grl_{\grg}$; let $\ad:U_{\grg}\to\End(\grg)$ be the $A$-homomorphism corresponding to the left $\grg$-module $\grg$. We refer to [CE, Ch. XIII] for the cohomology groups $H^i(\grg,\grv)$ of a left $\grg$-module $\grv$ (here $i\in\dbN$). We denote also by $\scrK_{\grg}:\grg\otimes_A\grg\to A$ the $A$-linear map defined by $\scrK_{\grg}:\grg\times \grg\to A$. Thus we have $\scrK_{\grg}\in (\grg\otimes_A \grg)^*=\grg^*\otimes_A\grg^*$. Let $\phi:\grg\arrowsim\grg^*$ be the $A$-linear isomorphism defined naturally by $\scrK_{\grg}$. It induces an $A$-linear isomorphism $\phi^{-1}\otimes\phi^{-1}:\grg^*\otimes_A \grg^*\arrowsim \grg\otimes_A \grg$. The image $\Omega$ of $\phi^{-1}\otimes\phi^{-1}(\scrK_{\grg})\in\grg\otimes_A \grg\subseteq T_{\grg}$ in $U_{\grg}$ is called the Casimir element of the adjoint representation $\ad:\grg\to\grg\grl_{\grg}$. 

\medskip\smallskip\noindent
{\bf 3.1. Lemma.} {\it For the Casimir element $\Omega\in U_{\grg}$ the following four properties hold:

\medskip
{\bf (a)} if the $A$-module $\grg$ is free and if $\{x_1,\ldots,x_m\}$ and $\{y_1,\ldots,y_m\}$ are two $A$-bases for $\grg$ such that for all $i,j\in\{1,\ldots,m\}$ we have $\scrK_{\grg}(x_i\otimes y_j)=\delta_{ij}$, then $\Omega$ is the image of the element $\sum_{i=1}^m x_i\otimes y_i$ of $T_{\grg}$ in $U_{\grg}$;

\smallskip
{\bf (b)} we have $\Omega\in Z(U_{\grg})$;

\smallskip
{\bf (c)} the Casimir element $\Omega$ is fixed by the group of Lie automorphisms of $\grg$ (i.e., if $\sigma:U_{\grg}\arrowsim U_{\grg}$ is the $A$-algebra automorphism induced by a Lie algebra automorphism $\sigma:\grg\arrowsim\grg$, then we have $\sigma(\Omega)=\Omega$); 

\smallskip
{\bf (d)} the Casimir element $\Omega$ acts identically on $\grg$ (i.e., $\ad(\Omega)=1_{\grg}$).}

\medskip
\proof
Parts (a) and (b) are proved in [B1, Ch. I, Sect. 3, Subsect. 7, Prop. 11]. Strictly speaking, loc. cit. is stated over a field but its proof applies over any commutative $\dbZ$-algebra. This is so as the essence of the proof of loc. cit. is [B1, Ch. I, Sect. 3, Subsect. 5, Example 2] which is worked out over any commutative $\dbZ$-algebra. In particular, [B1, Ch. I, Sect. 3, Subsect. 5, Example 3] can be easily stated over a commutative $\dbZ$-algebra (by involving a perfect invariant bilinear form over a commutative $\dbZ$-algebra instead of a non-degenerate invariant bilinear form over a field). We recall here that $\scrK_{\grg}$ is $\grg$-invariant i.e., for all $a,b,c\in\grg$ we have an identity $\scrK_{\grg}(\ad(a)(b),c)+\scrK_{\grg}(b,\ad(a)(c))=0$ (see [B1, Ch. I, Sect. 3, (13) and Prop. 8]) and this is the very essence of (b). 

To check (c) and (d), we can assume that the $A$-module $\grg$ is free. Let $\{x_1,\ldots,x_m\}$ and $\{y_1,\ldots,y_m\}$ be two $A$-bases for $\grg$ as in (a). Thus $\Omega$ is the image of $\sum_{i=1}^m x_i\otimes y_i$ in $U_{\grg}$. Therefore $\sigma(\Omega)$ is the image of $\sum_{i=1}^m \sigma(x_i)\otimes \sigma(y_i)$ in $U_{\grg}$. As for $i,j\in\{1,\ldots,m\}$ we have $\scrK_{\grg}(\sigma(x_i),\sigma(y_j))=\delta_{ij}$, from (a) we get that the image of $\sum_{i=1}^m \sigma(x_i)\otimes \sigma(y_i)\in T_{\grg}$ in $U_{\grg}$ is $\Omega$. Thus $\sigma(\Omega)=\Omega$.

We check (d). Let $z,w\in\grg$. We write $\ad(z)\circ \ad(w)(x_i)=\sum_{j=1}^m a_{ji}x_j$, with $a_{ji}$'s in $A$. Using the $\grg$-invariance of $\scrK_{\grg}$ we compute: 
$$\scrK_{\grg}(\ad(\Omega)(z),w)=\scrK_{\grg}(\sum_{i=1}^m\ad(x_i)\circ\ad(y_i)(z),w)=-\sum_{i=1}^m\scrK_{\grg}(\ad(y_i)(z),\ad(x_i)(w))$$
$$=\sum_{i=1}^m\scrK_{\grg}(\ad(z)(y_i),\ad(x_i)(w))=-\sum_{i=1}^m\scrK_{\grg}(y_i,\ad(z)\circ \ad(x_i)(w))$$
$$=\sum_{i=1}^m\scrK_{\grg}(y_i,\ad(z)\circ\ad(w)(x_i))=\sum_{i,j=1}^m a_{ji}\delta_{ji}=\sum_{i=1}^m a_{ii}=\scrK_{\grg}(z,w)$$
(the last equality due to the very definition of $\scrK_{\grg}$). This implies that for each $z\in\grg$, we have $\ad(\Omega)(z)-z\in\Ker(\scrK_{\grg})=0$. Thus $\ad(\Omega)(z)=z$ i.e., (d) holds.\endproof

\medskip\smallskip\noindent
{\bf 3.2. Fact.} {\it Let $i\in\dbN$. Let $\grv$ be a left $\grg$-module on which $\Omega$ acts identically. Then the cohomology group $H^i(\grg,\grv)$ is trivial.}

\medskip
\proof
We have an identity $H^i(\grg,\grv)=\text{Ext}^i_{U_{\grg}}(A,\grv)$ of $Z(U_{\grg})$-modules, cf. [CE, Ch. XIII, Sects. 2 and 8]. As $\Omega\in Z(U_{\grg})$ acts trivially on $A$ and identically on $\grv$, the group $\Omega \text{Ext}^i_{U_{\grg}}(A,\grv)$ is on one hand trivial and on the other hand it is $\text{Ext}^i_{U_{\grg}}(A,\grv)$. Thus $\text{Ext}^i_{U_{\grg}}(A,\grv)=0$. Therefore $H^i(\grg,\grv)=0$.\endproof

\medskip\smallskip\noindent
{\bf 3.3. Theorem.} {\it We recall that the Killing form $\scrK_{\grg}$ on $\grg$ is perfect. Then the group scheme $\text{Aut}(\grg)$ over $\Spec A$ of Lie automorphisms of $\grg$ is smooth and locally of finite presentation.}

\medskip
\proof
To check this, we can assume that the $A$-module $\grg$ is free. The group scheme $\text{Aut}(\grg)$ is a closed subgroup scheme of $\pmb{GL}_{\grg}$ defined by a finitely generated ideal of the ring of functions of $\pmb{GL}_{\grg}$. Thus $\text{Aut}(\grg)$ is of finite presentation. Thus to show that  $\text{Aut}(\grg)$ is smooth over $\Spec A$, it suffices to show that for each affine morphism $\Spec B\to\Spec A$ and for each ideal $\grj$ of $B$ such that $\grj^2=0$, the restriction map $\text{Aut}(\grg)(B)\to \text{Aut}(\grg)(B/\grj)$ is onto (cf. [BLR, Ch. 2, Sect. 2.2, Prop. 6]). Not to introduce extra notations by repeatedly tensoring with $B$ over $A$, we will assume that $B=A$. Thus $\grj$ is an ideal of $A$ and we have to show that the restriction map $\text{Aut}(\grg)(A)\to \text{Aut}(\grg)(A/\grj)$ is onto. 

Let $\bar\sigma:\grg/\grj\grg\arrowsim\grg/\grj\grg$ be a Lie automorphism. Let $\sigma_0:\grg\arrowsim\grg$ be an $A$-linear automorphism that lifts $\bar\sigma$. Let $\grj\grg_{\bar\sigma}$ be the left $\grg$-module which as an $A$-module is $\grj\grg$ and whose left $\grg$-module structure is defined as follows: if $x\in\grg$, then $x$ acts on $\grj\grg_{\bar\sigma}$ in the same way as $\ad(\bar\sigma(x))$ (equivalently, as $\ad(\sigma_0(x))$) acts on the $A$-module $\grj\grg=\grj\grg_{\bar\sigma}$; this makes sense as $\grj^2=0$. Let $\theta:\grg\times\grg\to \grj\grg_{\bar\sigma}$ be the alternating map defined by the rule:
$$\theta(x,y):=[\sigma_0(x),\sigma_0(y)]-\sigma_0([x,y])\,\;\;\;\forall\,x,y\in\grg.\leqno (1)$$
We check that $\theta$ is a $2$-cocycle i.e., for all $x,y,z\in\grg$ we have an identity
$$d\theta(x,y,z):=x(\theta(y,z))-y(\theta(x,z))+z(\theta(x,y))-\theta([x,y],z)+\theta([x,z],y)-\theta([y,z],x)=0.$$
Substituting (1) in the definition of $d\theta$, we get that the expression $d\theta(x,y,z)$ is a sum of 12 terms which can be divided into three groups as follows. The first group contains the three terms $-\sigma_0([[x,y],z])$, $\sigma_0([[x,z],y])$, and $-\sigma_0([[y,z],x])$; their sum is $0$ due to the Jacobi identity and the fact that $\sigma_0$ is an $A$-linear map. The second group contains the six terms $[\sigma_0(x),\sigma_0[y,z]]$, $-[\sigma_0(x),\sigma_0[y,z]]$, $[\sigma_0(y),\sigma_0[x,z]]$, $-[\sigma_0(y),\sigma_0[x,z]]$, $[\sigma_0(z),\sigma_0[x,y]]$, $-[\sigma_0(z),\sigma_0[x,y]]$; obviously their sum is $0$. The third group contains the three terms $[\sigma_0(x),[\sigma_0(y),\sigma_0(z)]]$, $-[\sigma_0(y),[\sigma_0(x),\sigma_0(z)]]$, and $[\sigma_0(z),[\sigma_0(x),\sigma_0(y)]]$; their sum is $0$ due to the Jacobi identity. Thus indeed $d\theta=0$.

As $\Omega$ (i.e., $\ad(\Omega)$) acts identically on $\grg$ (cf. Lemma 3.1 (d)), it also acts identically on $\grj\grg$. But $\Omega$ modulo $\grj$ is fixed by the Lie automorphism $\bar\sigma$ of $\grg/\grj\grg$, cf. Lemma 3.1 (c). Thus $\Omega$ also acts identically  on the left $\grg$-module $\grj\grg_{\bar\sigma}$. From this and the Fact 3.2 we get that $H^2(\grg,\grj\grg_{\bar\sigma})=0$. Thus $\theta$ is the coboundary of a $1$-cochain $\delta:\grg\to \grj\grg_{\bar\sigma}$ i.e., we have
$$\theta(x,y)=x(\delta(y))-y(\delta(x))-\delta([x,y])\;\;\;\;\forall x,\,y\in\grg.\leqno (2)$$
Let $\sigma:\grg\arrowsim\grg$ be the $A$-linear isomorphism defined by the rule $\sigma(x):=\sigma_0(x)-\delta(x)$; here $\delta(x)$ is an element of the $A$-module $\grj\grg=\grj\grg_{\bar\sigma}$. Due to formulas (1) and (2), we compute
$$\sigma([x,y])=\sigma_0([x,y])-\delta([x,y])=[\sigma_0(x),\sigma_0(y)]-\theta(x,y)-\delta([x,y])$$
$$=[\sigma_0(x),\sigma_0(y)]-x(\delta(y))+y(\delta(x))=
[\sigma_0(x),\sigma_0(y)]-\ad(\bar\sigma(x))(\delta(y))+\ad(\bar\sigma(y))(\delta(x))
$$
$$=[\sigma_0(x),\sigma_0(y)]-\ad(\sigma_0(x))(\delta(y))+\ad(\sigma_0)(y)(\delta(x))=[\sigma_0(x),\sigma_0(y)]-[\sigma_0(x),\delta(y)]+[\sigma_0(y),\delta(x)]$$
$$=[\sigma_0(x)-\delta(x),\sigma_0(y)-\delta(y)]-[\delta(x),\delta(y)]=[\sigma(x),\sigma(y)]-[\delta(x),\delta(y)]=[\sigma(x),\sigma(y)]$$
(the last identity as $\grj^2=0$). 
Thus $\sigma$ is a Lie automorphism of $\grg$ that lifts the Lie automorphism $\bar\sigma$ of $\grg/\grj\grg$. Thus the restriction map $\text{Aut}(\grg)(A)\to \text{Aut}(\grg)(A/\grj)$ is onto.\endproof

\medskip\smallskip\noindent
{\bf 3.5. Proof of the Theorem 1.2.}
The functor $\scrL_Y$ is faithful, cf. Lemma 2.5.1. Thus to prove Theorem 1.2 it suffices to show that $\scrL_Y$ is surjective on objects and that $\scrL_Y$ is fully faithful. To check this, as $\text{Adj-perf}_Y$ and $\text{Lie-perf}_Y$ are groupoids on sets and as $\scrL_Y$ is faithful, we can assume that $Y=\Spec A$ is affine. Thus to end the proof it suffices to check the following three properties: 

\medskip
{\bf (i)} if $\grg$ is an object of $\text{Lie-perf}_Y$ (identified with a Lie algebra over $A$), then there exists a unique open subgroup scheme $\text{Aut}(\grg)^0$ of $\text{Aut}(\grg)$ which is an adjoint group scheme over $Y$ and whose Lie algebra is the Lie subalgebra $\ad(\grg)$ of $\grg\grl_{\grg}$ (therefore $\grg=\ad(\grg)$ is the image through $\scrL_Y$ of the object $\text{Aut}(\grg)^0$ of $\text{Adj-perf}_Y$);

\smallskip
{\bf (ii)} the group scheme $\text{Aut}(\text{Aut}(\grg)^0)$ of automorphisms of $\text{Aut}(\grg)^0$ is $\text{Aut}(\grg)$ acting on $\text{Aut}(\grg)^0$ via inner conjugation (therefore $\text{Aut}(\grg)(A)=\text{Aut}(\text{Aut}(\grg)^0)(A)$);

\smallskip
{\bf (iii)} if $G$ and $H$ are two objects of $\text{Adj-perf}_Y$ such that $\Lie(G)=\Lie(H)$, then $G$ and $H$ are isomorphic. 

\medskip
To check the first two properties, we can assume that the $A$-module $\grg$ is free of rank $m\in\dbN^*$. Let $k$ be the residue field of an arbitrary point $y\in Y$. It is well known that the Lie algebra $\Lie(\text{Aut}(\grg)_k)$ is the Lie algebra of derivations of $\grg_k:=\grg\otimes_A k$. As this fact plays a key role in this paper, we include a proof of it. The tangent space of $\text{Aut}(\grg)_k$ at the identity element is identified with the set of automorphisms $a$ of $\grg_k\otimes_k k[\eps]/(\eps^2)$ which modulo $\bar\eps=\eps+(\eps^2)$ are the identity automorphism of $\grg_k$. We can write each such automorphism as $a=1_{\grg_k\otimes_k k[\eps]/(\eps^2)}+D_a\otimes\bar\eps$, where $D_a$ is a $k$-linear endomorphism of $\grg_k$. The condition that $a$ respects the Lie bracket (i.e., we have $a([u,v]\otimes 1)=[a(u\otimes 1),a(v\otimes 1)]$ for all $u,v\in\grg_k$) is equivalent to the condition that $D_a$ is a derivation of $\grg_k$. The association $a\mapsto D_a$ identifies the tangent space of $\text{Aut}(\grg)_k$ at the identity element with the $k$-vector space of derivations of $\grg_k$. Under this identification, the Lie bracket of $a$ with an automorphism $b$ of $\grg_k\otimes_k k[\eps_1]/(\eps_1^2)$ which modulo $\bar\eps_1=\eps_1+(\eps^2_1)$ is the identity automorphism of $\grg_k$, is the derivation of $\grg_k$ which corresponds to the automorphism $aba^{-1}b^{-1}=1_{\grg_k\otimes k[\eps\eps_1]/(\eps^2\eps_1^2)}+[D_a,D_b]\bar\eps\bar\eps_1$ of $\grg_k\otimes_k k[\eps\eps_1]/(\eps^2\eps_1^2)$ and thus is the Lie bracket $[D_a,D_b]$ ($\eps_1$ is used here instead of $\eps$ so that this last part makes sense). Therefore $\Lie(\text{Aut}(\grg)_k)$ is the Lie algebra of derivations of $\grg_k$.

As the Killing form $\scrK_{\grg_k}$ is perfect, as in [H1, Ch. II, Subsect. 5.3, Thm.] one argues that each derivation of $\grg_k$ is an inner derivation. Thus we have $\Lie(\text{Aut}(\grg)_k)=\ad(\grg)\otimes_A k$. As the group scheme $\text{Aut}(\grg)$ over $Y$ is smooth and locally of finite presentation (cf. Theorem 3.3), from [DG, Vol. I, Exp. VI${}_B$, Cor. 4.4] we get that there exists a unique open subgroup scheme $\text{Aut}(\grg)^0$ of $\text{Aut}(\grg)$ whose fibres are connected. The fibres of $\text{Aut}(\grg)^0$ are open-closed subgroups of the fibres of $\text{Aut}(\grg)$ and thus are affine.  

Let $N_k$ be a smooth, connected, unipotent, normal subgroup of $\text{Aut}(\grg)^0_k$. The Lie algebra $\Lie(N_k)$ is a nilpotent ideal of $\Lie(\text{Aut}(\grg)_k^0)=\ad(\grg)\otimes_A k$. Thus $\Lie(N_k)\subseteq \Ker(\scrK_{\Lie(\text{Aut}(\grg)_k^0)})=\Ker(\scrK_{\ad(\grg)\otimes_A k})$, cf. [B1, Ch. I, Sect. 4, Prop. 6 (b)]. As the Killing form $\scrK_{\ad(\grg)\otimes_A k}$ is perfect, we get $\Lie(N_k)=0$. Thus $N_k$ is the trivial subgroup of $\text{Aut}(\grg)_k^0$ and therefore the unipotent radical of $\text{Aut}(\grg)_k^0$ is trivial. Thus $\text{Aut}(\grg)_k^0$ is an affine, connected, smooth group over $\Spec k$ whose unipotent radical is trivial. Therefore $\text{Aut}(\grg)_k^0$ is a reductive group over $\Spec k$, cf. [Bo, Ch. IV, Subsect. 11.21]. As $\Lie(\text{Aut}(\grg)_k^0)=\ad(\grg)\otimes_A k$ has trivial center, the group $\text{Aut}(\grg)_k^0$ is semisimple. Thus the smooth group scheme $\text{Aut}(\grg)^0$ of finite presentation over $Y$ has semisimple fibres. Therefore $\text{Aut}(\grg)^0$ is a semisimple group scheme over $Y$, cf. [DG, Vol. II, Exp. XVI, Thm. 5.2 (ii)]. As $Z(\text{Aut}(\grg)^0)_k$ acts trivially on $\Lie(\text{Aut}(\grg)_k^0)=\ad(\grg)\otimes_A k$ and as $Z(\text{Aut}(\grg)^0)_k$ is a subgroup of $\text{Aut}(\grg)_k$, the group $Z(\text{Aut}(\grg)^0)_k$ is trivial. This implies that the finite, flat group scheme $Z(\text{Aut}(\grg)^0)$ is trivial and thus $\text{Aut}(\grg)^0$ is an adjoint group scheme. 

The Lie subalgebras $\Lie(\text{Aut}(\grg)^0)$ and $\ad(\grg)$ of $\grg\grl_{\grg}$ are free $A$-submodules of the Lie subalgebra $\grl$ of $\grg\grl_{\grg}$ formed by derivations of $\grg$. As for each point $y$ of $Y$ we have $\Lie(\text{Aut}(\grg)^0_k)=\ad(\grg)\otimes_A k=\grl\otimes_A k$, $\grl$ is locally generated by either $\Lie(\text{Aut}(\grg)^0)$ or $\ad(\grg)$. We easily get that we have identities $\Lie(\text{Aut}(\grg)^0)=\ad(\grg)=\grl$. 

The group scheme $\text{Aut}(\grg)$ acts via inner conjugation on $\text{Aut}(\grg)^0$. As $\Lie(\text{Aut}(\grg)^0)=\ad(\grg)$ and as $\text{Aut}(\grg)$ is a closed subgroup scheme of $\pmb{GL}_{\grg}$, the inner conjugation homomorphism $\text{Aut}(\grg)\to \Aut(\text{Aut}(\grg)^0)$ has trivial kernel. As $\Aut(\text{Aut}(\grg)^0)$ is a closed subgroup scheme of $\Aut(\Lie(\text{Aut}(\grg)^0))=\Aut(\ad(\grg))$ (cf.  Lemma 2.5.1), we can identify naturally $\Aut(\text{Aut}(\grg)^0)$ with a closed subgroup scheme of $\text{Aut}(\grg)$. From the last two sentences, we get that $\Aut(\text{Aut}(\grg)^0)=\text{Aut}(\grg)$. Thus both properties (i) and (ii) hold.

To check that the property (iii) holds, let $\grg=\Lie(G)=\Lie(H)$. It suffices to show that $G$ and $H$ are identified with $\text{Aut}(\grg)^0$. We will work only with $G$. The adjoint representation $G\to \pmb{GL}_{\grg}$ factors as composite closed embedding homomorphisms $G\to \text{Aut}(\grg)^0\to \text{Aut}(\grg)\to \pmb{GL}_{\grg}$ (cf. Lemma 2.5.1 and [DG, Vol. III, Exp. XXIV, Thm. 1.3]). We get a closed embedding homomorphism $G\to \text{Aut}(\grg)^0$ between adjoint group schemes that have the same Lie algebra $\grg$ (cf. also property (i)). By reasons of dimensions, the geometric fibers of the closed embedding homomorphism $G\to \text{Aut}(\grg)^0$ are isomorphisms and therefore $G\to \text{Aut}(\grg)^0$ is an isomorphism. Thus property (iii) holds as well.\endproof

\medskip
The next Proposition details on the range of applicability of the Theorem 1.2.

\medskip\smallskip\noindent
{\bf 3.6. Proposition.} {\it {\bf (a)} We recall that $k$ is a field. Let $H$ be a non-trivial semisimple group over $\Spec k$. Then the Killing form $\scrK_{\Lie(H)}$ is perfect if and only if the following two conditions hold:

\medskip
{\bf (i)} either $\char(k)=0$ or $\char(k)$ is an odd prime $p$ and $H^{\ad}$ has no simple factor of isotypic $A_{pn-1}$, $B_{pn+{{1-p}\over 2}}$, $C_{pn-1}$, or $D_{pn+1}$ Dynkin type (here $n\in\dbN^*$);

\smallskip
{\bf (ii)} if $\char(k)=3$ (resp. if $\char(k)=5$), then $H^{\ad}$ has no simple factor of isotypic $E_6$, $E_7$, $E_8$, $F_4$, $G_2$ (resp. of isotypic  $E_8$) Dynkin type.

\medskip
{\bf (b)} If $\scrK_{\Lie(H)}$ is perfect, then the central isogenies $H^{\sc}\to H\to H^{\ad}$ are \'etale; thus, by identifying tangent spaces at identity elements, we have $\Lie(H^{\sc})=\Lie(H)=\Lie(H^{\ad})$.}

\medskip
\proof
We can assume that $k=\bar k$ and that $\text{tr}.\text{deg}.(k)<\infty$. If $\char(k)=0$, then $\Lie(H)$ is a semisimple Lie algebra over $k$ and therefore the Proposition follows from [H1, Ch. II, Subsect. 5.1, Thm.]. Thus we can assume $\char(k)$ is a prime $p\in\dbN^*$. If the conditions (i) and (ii) hold, then $p$ does not divide the order of the finite group scheme $Z(H^{\sc})=\Ker(H^{\sc}\to H^{\ad})$ (see [B2, PLATES I to IX]) and therefore (a) implies (b). 

Let $W(k)$ be the ring of $p$-typical Witt vectors with coefficients in $k$. Let $H_{W(k)}$ be a semisimple group scheme over $\Spec W(k)$ that lifts $H$, cf. [DG,  Vol. III, Exp. XXIV, Prop. 1.21]. We have identities $\Lie(H^{\sc}_{W(k)})[{1\over p}]=\Lie(H_{W(k)})[{1\over p}]=\Lie(H^{\ad}_{W(k)})[{1\over p}]$. This implies that:

\medskip
{\bf (iii)} $\scrK_{\Lie(H_{W(k)})}$ is the composite of the natural $W(k)$-linear map $\Lie(H_{W(k)})\times \Lie(H_{W(k)})\to \Lie(H^{\ad}_{W(k)})\times \Lie(H^{\ad}_{W(k)})$ with $\scrK_{\Lie(H^{\ad}_{W(k)})}$;

\smallskip
{\bf (iv)} $\scrK_{\Lie(H^{\sc}_{W(k)})}$ is the composite of the natural $W(k)$-linear map $\Lie(H^{\sc}_{W(k)})\times \Lie(H^{\sc}_{W(k)})\to \Lie(H_{W(k)})\times \Lie(H_{W(k)})$ with $\scrK_{\Lie(H_{W(k)})}$. 

\medskip
We prove (a). We have $\Ker(\Lie(H)\to\Lie(H^{\ad}))\subseteq \Ker(\scrK_{\Lie(H)})$, cf. property (iii). If $\scrK_{\Lie(H)}$ is perfect, then $\Ker(\Lie(H)\to\Lie(H^{\ad}))=0$ and therefore $\Lie(H)=\Lie(H^{\ad})$. Thus to prove (a) we can assume that $H=H^{\ad}$ is adjoint. Even more, to prove (a) we can also assume that the adjoint group $H$ is simple; let $\flat$ be the Lie type of $H$. If $\flat$ is not of classical Lie type, then $\scrK_{\Lie(H)}$ is perfect if and only if either $p>5$ or $p=5$ and $\flat\neq E_8$ (cf. [H2, TABLE, p. 49]). Thus to prove (a), we can assume that $\flat$ is a classical Lie type.  We fix a morphism $\Spec \dbC\to\Spec W(k)$. 

Suppose that $\flat$ is either $A_n$ or $C_n$. By the standard trace form on $\Lie(H^{\sc})$ (resp. on $\Lie(H^{\sc}_{W(k)})$ or on $\Lie(H^{\sc}_{\dbC})$) we mean the trace form $\scrT$ (resp. $\scrT_{W(k)}$ or $\scrT_{\dbC}$) associated to the faithful representation of $H^{\sc}$ (resp. $H^{\sc}_{W(k)}$ or $H^{\sc}_{\dbC}$) of rank $n+1$ if $\flat=A_n$ and of rank $2n$ if $\flat=C_n$. We have $\scrK_{\Lie(H^{\sc}_{\dbC})}=2(n+1)\scrT_{\dbC}$, cf. [He, Ch. III, Sect. 8, (5) and (22)]. This identity implies that we also have $\scrK_{\Lie(H_{W(k)}^{\sc})}=2(n+1)\scrT_{W(k)}$ and thus $\scrK_{\Lie(H^{\sc})}=2(n+1)\scrT$. If $p$ does not divide $2(n+1)$, then $\Lie(H^{\sc})=\Lie(H)$ and it is well known that $\scrT$ is perfect; thus $\scrK_{\Lie(H^{\sc})}=\scrK_{\Lie(H)}=2(n+1)\scrT$ is perfect. Suppose that $p$ divides $2(n+1)$. This implies that $\scrK_{\Lie(H^{\sc})}$ is the trivial bilinear form on $\Lie(H^{\sc})$. From this and the property (iv) we get that the restriction of $\scrK_{\Lie(H)}$ to $\im(\Lie(H^{\sc})\to\Lie(H))$ is trivial. As $\dim_k(\Lie(H)/\im(\Lie(H^{\sc})\to\Lie(H)))=1$ and as $\dim_k(\Lie(H))\Ge 3$, we easily get that $\scrK_{\Lie(H)}$ is degenerate. 

Suppose that $\flat=B_n$ (resp. that $\flat=D_n$ with $n\Ge 4$). If $p>2$ we have $\Lie(H^{\sc})=\Lie(H)$. Moreover, using [He,  Ch. III, Sect. 8, (11) and (15)], as in the previous paragraph we argue that $\scrK_{\Lie(H)}$ is perfect if $p$ does not divide $2(2n-1)$ (resp. if $p$ does not divide $2(n-1)$) and is degenerate if $p$ divides $2n-1$ (resp. if $p$ divides $2(n-1)$). 

We are left to show that $\scrK_{\Lie(H)}$ is degenerate if $p=2$ and $\flat=B_n$. The group $H$ is (isomorphic to) the $\pmb{SO}$-group of the quadratic form $x_0^2+x_1x_{n+1}+\cdots+x_nx_{2n}$ on $W:=k^{2n+1}$. Let $\{e_{i,j}|i,j\in\{0,1,\ldots,n\}$ be the standard $k$-basis for $\gl_W$. The direct sum $\grn_n:=\oplus_{i=1}^{2n} ke_{0,i}$ is a nilpotent ideal of $\Lie(H)$, cf. [Bo, Ch. V, Subsect. 23.6]. Thus $\grn_n\subseteq \Ker(\scrK_{\Lie(H)})$, cf. [B1, Ch. I, Sect. 4, Prop. 6 (b)] applied to the adjoint representation of $\Lie(H)$. Therefore $\scrK_{\Lie(H)}$ is degenerate. 

We conclude that $\scrK_{\Lie(H)}$ is perfect if and only if both conditions (i) and (ii) hold. Therefore (a) (and thus also (b)) holds.\endproof

\medskip\smallskip\noindent
{\bf 3.7. Remarks.} Let $A$ and $\grg$ be as in the beginning of this Section.

\smallskip
{\bf (a)} Let $p\in\dbN^*$ be a prime. Suppose that $A$ is an algebraically closed field of characteristic $p$. Let $G$ be an adjoint group over $\Spec A$ such that $\grg=\Lie(G)$, cf. Theorem 1.2. We have $p\neq 2$, cf. Proposition 3.6. Let $G_{\dbZ}$ be the unique (up to isomorphism) split, adjoint group scheme over $\Spec \dbZ$ such that $G$ is the pull back of $G_{\dbZ}$ to $\Spec A$, cf. [DG,  Vol. III, Exp. XXV, Cor. 1.3]. We have $\grg=\Lie(G_{\dbZ})\otimes_{\dbZ} A$ i.e., $\grg$ has a canonical model $\Lie(G_{\dbZ})$ over $\dbZ$. For $p>7$, this result was obtained in [C, Sect. 5, Thm.]. For $p>3$, this result was obtained by Seligman, Mills, Block, and Zassenhaus (see [MS], [Mi], [BZ], and [S, II. 10]). For $p=3$, this result was obtained in [Br, Thm. 4.1]. It seems to us that the fact that $p\neq 2$ (i.e., that all Killing forms of finite dimensional Lie algebras over fields of characteristic $2$, are degenerate) is new. 

\smallskip
{\bf (b)} Let $B\twoheadrightarrow A$ be an epimorphism of commutative $\dbZ$-algebras whose kernel $\grj$ is a nilpotent ideal. Then $\grg$ has, up to isomorphisms, a unique lift to a Lie algebra over $B$ which as a $B$-module is projective and finitely generated. One can prove this statement using cohomological methods as in the proof of Theorem 3.3. The statement also follows from the Theorem 1.2 and the fact that $\Aut(\grg)^0$ has, up to isomorphisms, a unique lift to an adjoint group scheme over $\Spec B$ (this can be easily checked at the level of torsors of adjoint group schemes; see [DG, Vol. III, Exp. XXIV, Cors. 1.17 and 1.18]). 

\medskip\smallskip\noindent
{\bf 3.8. Corollary.} {\it Let $\text{Sc-perf}_Y$ be the category whose objects are simply connected semisimple group schemes over $Y$ with the property that their Lie algebras $\scrO_Y$-modules have perfect Killing forms and whose morphisms are isomorphisms of group schemes. Then the functor $\scrL^{\text{sc}}_Y:\text{Sc-perf}_Y\to\text{Lie-perf}_Y$ which associates to a morphism $f:G\arrowsim H$ of $\text{Sc-perf}_Y$ the morphism $df:Lie(G)\arrowsim Lie(H)$ of $\text{Lie-perf}_Y$ which is the differential of $f$, is an equivalence of categories.}

\medskip
\proof
The functor $\scrL^{\text{sc}}_Y$ is the composite of the canonical (`division by the centers') functor $\scrZ_Y:\text{Sc-perf}_Y\to \text{Ad-perf}_Y$ with $\scrL_Y$; the functor $\scrZ_Y$ makes sense (cf. Lemma 3.6 (b)) and it is an equivalence of categories. Thus the Corollary follows from the Theorem 1.2.\endproof

\medskip\smallskip\noindent
{\bf 3.9. Corollary.} {\it The category $\text{Lie-perf}_Y$ has a non-zero object if and only if $Y$ is a non-empty $\Spec \dbZ[{1\over 2}]$-scheme.}

\medskip
\proof
The if part is implied by the fact that an $\grs\grl_2$ Lie algebra $\scrO_Y$-module has perfect Killing form. The only if part follows from the relation $p\neq 2$ of the Remark 3.7 (a).\endproof

\bigskip
\noindent
{\boldsectionfont 4. Proof of the Theorem 1.4}
\bigskip 

In this Section we prove Theorem 1.4. See Subsections 4.1 and 4.2 for the proofs of Theorems 1.4 (a) and (b) (respectively). In Remarks 4.3 we point out that the hypotheses of the Theorem 1.4 are indeed needed in general. We will use the notations listed in Section 1.

\medskip\smallskip\noindent
{\bf 4.1. Proof of Theorem 1.4 (a).} To prove Theorem 1.4  we can assume $Y$ is also integral. Let $K:=K_Y$; it is a field. If $H$ is a reductive group scheme over $Y$, then we have $\scrD(H)=\scrD(H_U)$ and thus the uniqueness parts of the Theorem 1.4 follow from Proposition 2.2 (b). Let $\grl$ be the Lie algebra $\scrO_Y$-module which extends $\Lie(G_U)$.

We prove Theorem 1.4 (a). Due to the uniqueness part, to prove Theorem 1.4 (a) we can assume $Y=\Spec A$ is also local and strictly henselian. Let $\grg:=\grl(Y)$ be the Lie algebra over $A$ of global sections of $\grl$.

As $U$ is connected, based on [DG, Vol. III, Exp. XXII, Prop. 2.8] we can speak about the split, adjoint group scheme $S$ over $Y$ of the same Lie type as all geometric fibres of $G_U$. Let $\grs:=\Lie(S)$. Let $\Aut(S)$ be the group scheme over $Y$ of automorphisms of $S$. We have a short exact sequence $1\to S\to \text{Aut}(S)\to C\to 1$, where $C$ is a finite, \'etale, constant  group scheme over $Y$ (cf. [DG, Vol. III, Exp. XXIV, Thm. 1.3]). Let $\gamma\in H^1(U,\text{Aut}(S)_U)$ be the class that defines the form $G_U$ of $S_U$. 

We recall that $\pmb{GL}_{\grg}$ and $\pmb{GL}_{\grs}$ are the reductive group schemes over $Y$ of linear automorphisms of $\grg$ and $\grs$ (respectively). The adjoint representations define closed embedding homomorphisms $j_U:G_U\hookrightarrow \pmb{GL}_{\grg,U}$ and $i:S\hookrightarrow \pmb{GL}_{\grs}$ and moreover $i$ extends naturally to a closed embedding homomorphism $\Aut(S)\hookrightarrow \pmb{GL}_{\grs}$, cf. Lemma 2.5.1. Let $\delta\in H^1(U,(\pmb{GL}_{\grs,U}))$ be the image of $\gamma$ via the homomorphism $\Aut(S)_U\hookrightarrow \pmb{GL}_{\grs,U}$. 

We recall that the quotient sheaf for the faithfully flat topology of $Y$ of the action of $S$ on $\pmb{GL}_{\grs}$ via right translations, is representable by an $Y$-scheme $\pmb{GL}_{\grs}/S$ that is affine and that makes $\pmb{GL}_{\grs}$ to be a right torsor of $S$ over $\pmb{GL}_{\grs}/S$ (cf. [CTS, Cor. 6.12]). Thus $\pmb{GL}_{\grs}/S$ is a smooth, affine $Y$-scheme. The finite, \'etale, constant group scheme $C$ acts naturally (from the right) on $\pmb{GL}_{\grs}/S$ and this action is free (cf. Lemma 2.5.1). From [DG, Vol I, Exp. V, Thm. 4.1] we get that the quotient $Y$-scheme $(\pmb{GL}_{\grs}/S)/C$ is affine and that the quotient epimorphism $\pmb{GL}_{\grs}/S\twoheadrightarrow (\pmb{GL}_{\grs}/S)/C$ is a finite \'etale cover. Thus $(\pmb{GL}_{\grs}/S)/C$ is a smooth, affine scheme over $Y$ that represents the quotient sheaf for the faithfully flat topology of $Y$ of the action of $\Aut(S)$ on $\pmb{GL}_{\grs}$ via right translations. From constructions we get that $\pmb{GL}_{\grs}$ is a right torsor of $\Aut(S)$ over $\pmb{GL}_{\grs}/\Aut(S):=(\pmb{GL}_{\grs}/S)/C$.

The twist of $i_U$ via the class $\gamma$ is $j_U$. This implies that the class $\delta$ defines the torsor that parametrizes isomorphisms between the pull backs to $U$ of the vector group schemes over $Y$ defined by $\grs$ and $\grg$. Therefore, as the $A$-modules $\grs$ and $\grg$ are isomorphic (being free of equal ranks), the class $\delta$ is trivial. Thus $\gamma$ is the coboundary of a class in $H^0(U,\pmb{GL}_{\grs,U}/\Aut(S)_U)$. But $H^0(U,\pmb{GL}_{\grs,U}/\Aut(S)_U)=H^0(Y,\pmb{GL}_{\grs}/\Aut(S))$ (cf. Proposition 2.2 (b)) and therefore $\gamma$ is the restriction of a class in $H^1(Y,\Aut(S))$. As $Y$ is strictly henselian, each class in $H^1(Y,\Aut(S))$ is trivial. Thus $\gamma$ is the trivial class. Therefore the group schemes $G_U$ and $S_U$ are isomorphic. Thus $G_U$ extends to an adjoint group scheme $G$ over $Y$ isomorphic to $S$. This ends the proof of Theorem 1.4 (a). \endproof

\medskip\smallskip\noindent
{\bf 4.2. Proof of Theorem 1.4 (b).} Let $\eta:\bar K\to U$ be the geometric point of $U$ which is the composite of the natural morphisms $\Spec \bar K\to\Spec K$ and $\Spec K\to U$. We denote also by $\eta:\bar K\to Y$ the resulting geometric point of $Y$. As $Y$ (resp. $U$) is normal and locally noetherian, from [DG, Vol. II, Exp. X, Thms. 5.16 and 7.1] we get that there exists an antiequivalence of categories between the category of tori over $Y$ (resp. over $U$) and the category of continuous $\pi_1(Y,\eta)$-representations (resp. continuous $\pi_1(U,\eta)$-representations) on free $\dbZ$-modules of finite rank. As the pair $(Y,Y\setminus U)$ is quasi-pure, we have a canonical identification $\pi_1(U,\eta)=\pi(Y,\eta)$. From the last two sentences we get that there exists a unique torus $H^{\ab}$ over $Y$ which extends $H_U^{\ab}$.

Let $H^{\ad}$ be the adjoint group scheme over $Y$ that extends $H^{\ad}_U$, cf. Theorem 1.4 (a). Let $F\to H^{\ad}\times_Y H^{\ab}$ be the central isogeny over $Y$ that extends the central isogeny $H_K\to H_K^{\ad}\times_{\Spec K} H_K^{\ab}$, cf. Lemma 2.3.1 (a). Both $F_U$ and $H_U$ are the normalization of $H_U^{\ab}\times_U H_U^{\ad}$ in $H_K$, cf. Lemma 2.3.1 (b). Thus $H_U=F_U$ extends uniquely to a reductive group scheme $H:=F$ over $Y$ (cf. the first paragraph of Subsection 4.1 for the uniqueness part). This ends the proof of Theorem 1.4 (b) and thus also of the Theorem 1.4.\endproof
 
\medskip\smallskip\noindent
{\bf 4.3. Remarks.} {\bf (a)} Let $Y_1\to Y$ be a finite, non-\'etale morphism between normal, noetherian, integral $\Spec \dbZ_{(2)}$-schemes such that there exists an open subscheme $U$ of $Y$ with the properties that: (i) $Y\setminus U$ has codimension in $Y$ at least $2$, and (ii) $Y_1\times_Y U\to U$ is a Galois cover of degree $2$. Let $H_U$ be the rank $1$ non-split torus over $U$ that splits over $Y_1\times_Y U$. Then $H_U$ does not extend to a smooth, affine group scheme over $Y$. If moreover $Y=\Spec A$ is an affine $\Spec \dbF_2$-scheme, then we have $Lie(H_U)(U)=A$ and therefore $Lie(H_U)$ extends to a Lie algebra $\scrO_Y$-module which as an $\scrO_Y$-module is free. Thus the quasi-pure part of the hypotheses of Theorem 1.4 (b) is needed in general.

{\bf (b)} Suppose that $Y=\Spec A$ is local, strictly henselian, regular, and of dimension $n\Ge 3$. Let $K:=K_Y$. Let $d\in\dbN^*$ be such that there exists an $A$-submodule $M$ of $K^d$ that contains $A^d$, that is of finite type, that is not free, and that satisfies the identity $M=\cap_{V\in\scrD(Y)} M\otimes_A V$ (inside $M\otimes_A K$). A typical example (communicated to us by Serre): $d=n-1$ and $M\arrowsim\text{Coker}(f)$, where the $A$-linear map $f:A\to A^n$ takes $1$ to an $n$-tuple $(x_1,\ldots,x_n)\in A^n$ of regular parameters of $A$. 

Let $\scrF$ be the coherent $\scrO_Y$-module defined by $M$. Let $U$ be an open subscheme of $Y$ such that $Y\setminus U$ has codimension in $Y$ at least $2$ and the restriction $\scrF_U$ of $\scrF$ to $U$ is a locally free $\scrO_U$-module. Let $H_U$ be the reductive group scheme over $U$ of linear automorphisms of $\scrF_U$. We recall the reason why the assumption that $H_U$ extends to a reductive group scheme $H$ over $Y$ leads to a contradiction. The group scheme $H$ is isomorphic to $\pmb{GL}_{d,A}$ (as $A$ is strictly henselian) and therefore there exists a free $A$-submodule $L$ of $K^d$ of rank $d$ such that we can identify $H=\pmb{GL}_L$.  As $A$ is a unique factorization domain (being local and regular), it is easy to see that there exists an element $f\in K$ such that the identity $M\otimes_A V=fL\otimes_A V$ holds for each $V\in\scrD(Y)$. This implies that $M=fL$. Thus $M$ is a free $A$-module. Contradiction. 

As $H_U$ does not extend to a reductive group scheme over $Y$ and as the pair $(Y,Y\setminus U)$ is quasi-pure, from Subsection 4.2 we get that $H_U^{\ad}$ also does not extend to an adjoint group scheme over $Y$. Thus the Lie part of the hypotheses of Theorem 1.4 (a) is needed in general.

\bigskip
\noindent
{\boldsectionfont 5. Extending homomorphisms via schematic closures}
\bigskip 

In this Section we prove two results on extending homomorphisms of reductive group schemes via taking (normalizations of) schematic closures. The first one complements Theorem 1.4 (b) and Proposition 2.5.2 (see Proposition 5.1) and the second one refines [V1, Lemma 3.1.6] (see Proposition 5.2). 

\medskip\smallskip\noindent
{\bf 5.1. Proposition.} {\it Let $Y$ be a normal, noetherian, integral scheme. Let $K:=K_Y$. Let $U$ be an open subscheme of $Y$ such that the codimension of $Y\setminus U$ in $Y$ is at least $2$. Let $H_U$ be a reductive group scheme over $U$ and let $G$ be a reductive group scheme over $Y$. We assume we have a finite homomorphism $\rho_U:H_U\to G_U$ whose generic fiber over $\Spec K$ is a closed embedding. We assume that one of the following two properties holds:

\medskip
{\bf (i)} $H_U$ extends to a reductive group scheme $H$ over $Y$; 

\medskip
{\bf (ii)} $Y=\Spec R$ is a local regular scheme of dimension $2$ (thus $U$ is the complement in $Y$ of the closed point of $Y$).

\medskip
Then the following three properties hold:

\medskip
{\bf (a)} There exists a unique reductive group scheme $H$ over $Y$ which extends $H_U$.

\medskip
{\bf (b)} The homomorphism $\rho_U$ extends uniquely to a finite homomorphism $\rho:H\to G$ between reductive group schemes over $Y$. 

\medskip
{\bf (c)} If there exists a point of $Y\setminus U$ of characteristic $2$, we assume that $H_K$ has no normal subgroup that is adjoint of isotypic $B_n$ Dynkin type for some $n\in\dbN^*$. Then $\rho:H\to G$ is a closed embedding.}

\medskip
\proof
If the property (i) holds, then the uniqueness of $H_U$ follows from Proposition 2.2 (b). Thus to prove (a) we can assume that the property (ii) holds.  As (ii) holds, the pair $(Y,Y\setminus U)$ is quasi-pure (see Section 1) and the Lie algebra $\scrO_U$-module $Lie(H_U)$  extends to a Lie algebra $\scrO_Y$-module which is a free $\scrO_Y$-module (cf. Proposition 2.2 (c) and the fact that $Y$ is local). Thus the hypotheses of Theorem 1.4 (b) hold and therefore from Theorem 1.4 (b) we get that there exists a unique reductive group scheme $H$ over $Y$ that extends $H_U$. Thus (a) holds.

To prove (b) and (c) we can assume that $Y=\Spec R$ is an affine scheme. We write $H=\Spec R_H$ and $G=\Spec R_G$. As $\scrD(H)=\scrD(H_U)$ and $\scrD(G)=\scrD(G_U)$, from Proposition 2.2 (a) we get that $R_H$ and $R_G$ are the $R$-algebras of global functions of $H_U$ and $G_U$ (respectively). Let $R_G\to R_H$ be the $R$-homomorphism defined by $\rho_U$ and let $\rho:H\to G$ be the morphism of $Y$-schemes it defines. The morphism $\rho$ is a homomorphism as it is so generically. To check that $\rho$ is finite, we can assume that $R$ is complete. Thus $R_H$ and $R_G$ are excellent rings, cf. [M, Sect. 34]. Therefore the normalization $H^\prime=\Spec R_{H^\prime}$ of the schematic closure of $H_K$ in $G$ is a finite, normal $G$-scheme. 

The identity components of the reduced geometric fibres of $\rho$ are trivial groups, cf. Proposition 2.5.2 (a) or (b). Thus $\rho$ is a quasi-finite morphism. From Zariski Main Theorem (see [G1, Thm. (8.12.6)]) we get that $H$ is an open subscheme of $H^\prime$. But from Proposition 2.5.2 (b) we get that the morphism $H\to H^\prime$ satisfies the valuative criterion of properness with respect to discrete valuation rings which contain $R$. As each local ring of $H^\prime$ is dominated by such a discrete valuation ring, we get that the morphism $H\to H^\prime$ is surjective. Therefore the open, surjective morphism $H\to H^\prime$ is an isomorphism. Thus $\rho$ is finite i.e., (b) holds. 

We prove (c). The pull back of the homomorphism $\rho:H\to G$ via each dominant morphism $\Spec V\to Y$, with $V$ a discrete valuation ring, is a closed embedding (cf. Proposition 2.5.2 (c)). This implies that the fibres of $\rho$ are closed embeddings. Thus the homomorphism $\rho$ is a closed embedding, cf. Theorem 2.5.\endproof

\medskip
We have the following refinement of [V1, Lemma 3.1.6].

\medskip\smallskip\noindent
{\bf 5.2. Proposition.} {\it Let $G$ be a reductive group scheme over a reduced, affine scheme $Y=\Spec A$. Let $K$ be a localization of $A$. Let $s\in\dbN^*$. For $j\in\{1,\ldots,s\}$ let $G_{j,K}$ be a reductive, closed subgroup scheme of $G_K$. We assume that the group subschemes $G_{j,K}$'s commute among themselves and that one of the following two conditions holds:

\medskip
{\bf (i)} either the direct sum $\oplus_{j=1}^s \Lie(G_{j,K})$ is a Lie subalgebra of $\Lie(G_K)$, or 

\smallskip
{\bf (ii)} $s=2$, $G_{1,K}$ is a torus, and $G_{2,K}$ is  a semisimple group scheme. 

\medskip
Then the closed subgroup scheme $G_{0,K}$ of $G_K$ generated by $G_{j,K}$'s exists and is reductive. Moreover, we have:

\medskip
{\bf (a)} If the condition (i) holds, then $\Lie(G_{0,K})=\oplus_{j=1}^s \Lie(G_{j,K})$. 

\smallskip
{\bf (b)} We assume that for each $j\in\{1,\ldots,s\}$ the schematic closure $G_j$ of $G_{j,K}$ in $G$ is a reductive group scheme over $Y$. Then the schematic closure $G_0$ of $G_{0,K}$ in $G$ is a reductive, closed subgroup scheme of $G$.}

\medskip
\proof
Let $\Lambda$ be the category whose objects $Ob(\Lambda)$ are finite subsets of $K$ and whose morphisms are the inclusions of subsets. For $\alpha\in Ob(\Lambda)$, let $K_{\alpha}$ be the $\dbZ$-subalgebra of $K$ generated by $\alpha$ and let $A_{\alpha}:=A\cap K_{\alpha}$. We have $K=\text{ind.}\,\text{lim.}_{\alpha\in Ob(\Lambda)} K_{\alpha}$ and $A=\text{ind.}\,\text{lim.}_{\alpha\in Ob(\Lambda)} A_{\alpha}$. The reductive group schemes $G_{j,K}$ are of finite presentation. Based on this and [G1, Thms. (8.8.2) and (8.10.5)], one gets that there there exists $\beta\in Ob(\Lambda)$ such that each $G_{j,K}$ is the pull back of a closed subgroup scheme $G_{j,K_{\beta}}$ of $G_{K_{\beta}}$. For $\alpha\supseteq\beta$, the set $C(\alpha)$ of points of $\Spec K_{\alpha}$ with the property that the fibres over them of all morphisms $G_{j,K_{\alpha}}\to\Spec K_{\alpha}$ are (geometrically) connected, is a constructible set (cf. [G1, Thm. (9.7.7)]). We have $\text{proj.}\,\text{lim.}_{\alpha\in Ob(\Lambda)} C(\alpha)=\Spec K$. From this and [G1, Thm. (8.5.2)], we get that there exists $\beta_1\in Ob(\Lambda)$ such that $\beta_1\supseteq\beta$ and $C(\beta_1)=\Spec K_{\beta_1}$. Thus by replacing $\beta$ with $\beta_1$, we can assume that the fibres of all morphisms $G_{j,K_{\beta}}\to\Spec K_{\beta}$ are connected. A similar argument shows that, by enlarging $\beta$, we can assume that all morphisms $G_{j,K_{\beta}}\to\Spec K_{\beta}$ are smooth and their fibres are reductive groups (the role of [G1, Thm. (9.7.7)] being replaced by [G1, Prop. (9.9.5)] applied to the $\scrO_{G_{j,K_{\alpha}}}$-module $Lie(G_{j,K_{\alpha}})$ and respectively by [DG, Vol. III, Exp. XIX, Cor. 2.6]). Thus each $G_{j,K_{\beta}}$ is a reductive closed subgroup scheme of $G_{K_{\beta}}$. The smooth group schemes $G_{j,K_{\beta}}$'s commute among themselves as this is so after pull back through the dominant morphism $\Spec K\to\Spec K_{\beta}$. By enlarging $\beta$, we can also assume that either condition (i) or condition (ii) holds for the $G_{j,K_{\beta}}$'s and that $K_{\beta}$ is a localization of $A_{\beta}$. By replacing $A$ with the local ring of $\Spec A_{\beta}$ dominated by $A$, to prove the Proposition we can assume that $A$ is a localization of a reduced, finitely generated $\dbZ$-algebra.

Using induction on $s\in\dbN^*$, it suffices to prove the Proposition for $s=2$. Moreover, we can assume that $K=K_Y$. For the sake of flexibility, in what follows we will only assume that $A$ is a reduced, noetherian $\dbZ$-algebra; thus $K$ is a finite product of fields. As all the statements of the Proposition are local for the \'etale topology of $Y$, it suffices to prove the Proposition under the extra assumption that $G_1$ and $G_2$ are split (cf. Proposition 2.3). Let $C_K:=G_{1,K}\cap G_{2,K}$. It is a closed subgroup scheme of $G_{j,K}$ that commutes with $G_{j,K}$, $j\in\{1,2\}$. The Lie algebra $\Lie(C_K)$ is included in $\Lie(G_{1,K})\cap\Lie(G_{2,K})$ and therefore it is trivial if the condition (i) holds. Thus if the condition (i) holds, then $C_K$ is a finite, \'etale, closed subgroup scheme of $Z(G_{j,K})$. If the condition (ii) holds, then $C_K$ is a closed subgroup scheme of both $G_{1,K}=Z(G_{1,K})$ and $Z(G_{2,K})$ and thus (as $K$ is a finite product of fields) it is a finite group scheme of multiplicative type. 

Let $C$ be the schematic closure of $C_K$ in $G$. Let $T_j$ be a maximal torus of $G_j$. We have $C_K\leqslant Z(G_{1,K})\cap Z(G_{2,K})\leqslant T_{1,K}\cap T_{2,K}\leqslant G_{1,K}\cap G_{2,K}=C_K$ and thus $C_K=T_{1,K}\cap T_{2,K}$. Let $T_1\times_Y T_2\to G$ be the product homomorphism. The kernel $\grK$ of this product homomorphism is a group scheme over $Y$ of multiplicative type (cf. Lemma 2.3.2 (a)) isomorphic to $T_1\cap T_2$. But $\grK_K\arrowsim C_K$ is a finite group scheme over $\Spec K$ and therefore $\grK$ is a finite, flat group scheme over $Y$ of multiplicative type (cf. Lemma 2.3.2 (b)). Thus $T_1\cap T_2$ is a finite, flat group scheme over $Y$. From this, the identity $C_K=(T_1\cap T_2)_K$, and the definition of $C$ we get that $C=T_1\cap T_2$. We conclude that $C$ is a finite, flat group scheme over $Y$ of multiplicative type contained in the center of both $G_1$ and $G_2$. We embed $C$ in $G_1\times_Y G_2$ via the natural embedding $C\hookrightarrow G_1$ and via the composite of the inverse isomorphism $C\arrowsim C$ with the natural embedding $C\hookrightarrow G_2$. Let $G_{1,2}:=(G_1\times_Y G_2)/C$; it is a reductive group scheme over $Y$. We have a natural product homomorphism $q:G_{1,2}\to G$ whose pull back to $\Spec K$ can be identified with the closed embedding homomorphism $G_{0,K}\hookrightarrow G_K$. Therefore $G_{0,K}$ is a reductive group scheme over $\Spec K$. Moreover, if  the condition (i) holds, then as $C_K$ is \'etale we have natural identities $\Lie(G_{1,K})\oplus\Lie(G_{2,K})=\Lie(G_{1,2,K})=\Lie(G_{0,K})$. Thus (a) holds. If $q$ is a closed embedding, then $q$ induces an isomorphism $G_{1,2}\arrowsim G_0$ and therefore $G_0$ is a reductive, closed subgroup scheme of $G$. Thus to end the proof of (b), we only have to show that the homomorphism $q$ is a closed embedding.

To check that $q$ is a closed embedding, it suffices to check that the fibres of $q$ are closed embeddings (cf. Theorem 2.5). For this we can assume that $A$ is a complete discrete valuation ring which has an algebraically closed residue field $k$; this implies that $G_0$ is a flat, closed subgroup scheme of $G$. Let $\grn:=\Lie(\Ker(q_k))$. From Proposition 2.5.2 (a) and Lemma 2.4 we get that: either (iii) $\grn=0$ or (iv) $\char(k)=2$ and there exists a normal subgroup $F_k$ of $G_{1,2,k}$ which is isomorphic to $\pmb{SO}_{2n+1,k}$ for some $n\in\dbN^*$ and for which we have $\Lie(F_k)\cap\grn\neq 0$. We show that the assumption that the condition (iv) holds leads to a contradiction. Let $F$ be a normal, closed subgroup scheme of $G_{1,2}$ that lifts $F_k$ and that is isomorphic to $\pmb{SO}_{2n+1,A}$ (cf. last paragraph of the proof of Proposition 2.5.2 (c)). Let $j_0\in\{1,2\}$ be such that $F\vartriangleleft G_{j_0}\vartriangleleft G_{1,2}$ (if the condition (ii) holds, then $j_0=2$). As $G_{j_0}$ is a closed subgroup scheme of $G$, we have $\Lie(G_{j_0,k})\cap\grn=0$ and therefore also $\Lie(F_k)\cap\grn=0$. Contradiction. Thus the condition (iv) does not hold and therefore the condition (iii) holds. Thus $\Ker(q_k)$ has a trivial Lie algebra and therefore it is a finite, \'etale, normal subgroup of $G_{1,2,k}$. Thus $\Ker(q_k)$ is a subgroup of $Z(G_{1,2,k})$ and therefore also of each maximal torus of $G_{1,2,k}$. From this and Proposition 2.5.2 (a) we get that $\Ker(q_k)$ is trivial. Therefore $q_k$ is a closed embedding. Thus $q$ is a closed embedding.\endproof 

\medskip\noindent
{\bf Acknowledgments.} We would like to thank University of Utah, University of Arizona, Binghamton University, and I. H. E. S., Bures-sur-Yvette for good working conditions. We would like to thank J.-P. Serre and A. Langer for pointing out the reference [CTS] and O. Gabber for an e-mail that independently outlined parts of [CTS]; these led to a shorter paper and to stronger forms of Theorems 1.2 and 1.4 and of Proposition 5.1, for those reductive group schemes whose adjoints have geometric fibres that have simple factors of exceptional Lie types. We would like to thank J.-L. Colliot-Th\'el\`ene for some comments and J. E. Humphreys for pointing out the references quoted in Remark 3.7 (a). We would like to thank O. Gabber also for many other valuable comments and suggestions; Section 3 is mainly based on these suggestions. This research was partially supported by the NSF grant DMS \#0900967.

\bigskip
\references{37}
{\nspace{

\Ref[B1]
N. Bourbaki,
\sl Lie groups and Lie algebras, Chapters {\bf 1--3},
\rm Elements of Mathematics (Berlin), Springer-Verlag, Berlin, 1989.

\Ref[B2]
N. Bourbaki,
\sl Lie groups and Lie algebras, Chapters {\bf 4--6},
\rm Elements of Mathematics (Berlin), Springer-Verlag, Berlin, 2002.

\Ref[BLR]
S. Bosch, W. L\"utkebohmert, and M. Raynaud,
\sl N\'eron models,
\rm Ergebnisse der Mathematik und ihrer Grenzgebiete (3), Band {\bf 21}, Springer-Verlag, Berlin, 1990.

\Ref[BZ]
R. Block and H. Zassenhaus,
\sl The Lie algebras with a non-degenerate trace form,
\rm Illinois J. Math. {\bf 8} (1964),  543--549. 

\Ref[Bo]
A. Borel,
\sl Linear algebraic groups. Second edition,
\rm Grad. Texts in Math., Vol. {\bf 126}, Springer-Verlag, New York, 1991.

\Ref[Br]
G. Brown,
\sl Lie algebras of characteristic three with nondegenerate Killing form,
\rm Trans. Amer. Math. Soc. {\bf 137} (1969),  259--268.

\Ref[C]
C. W. Curtis,
\sl Modular Lie algebras II,
\rm Trans. Amer. Math. Soc., Vol. {\bf 86} (1957), no. 1,  91--108. 

\Ref[CE]
H. Cartan and S. Eilenberg,
\sl Homological algebra,
\rm Reprint of the 1956 original, Princeton Landmarks in Mathematics, Princeton Univ. Press, Princeton, NJ, 1999. 

\Ref[CTS]
J.-L. Colliot-Th\'el\`ene et J.-J. Sansuc,
\sl Fibr\'es quadratiques et composantes connexes r\'eelles,
\rm Math. Ann. {\bf 244} (1979), no. 2,  105--134. 

\Ref[DG]
M. Demazure, A. Grothendieck, et al.,
\sl Sch\'emas en groupes. Vols. {\bf I--III},
\rm S\'eminaire de G\'eom\'etrie Alg\'ebrique du Bois Marie 1962/64 (SGA 3), Lecture Notes in Math., Vols. {\bf 151--153}, Springer-Verlag, Berlin-New York, 1970.

\Ref[dJ] 
J. de Jong, 
\sl Homomorphisms of Barsotti--Tate groups and crystals in positive characteristic,
\rm Invent. Math. {\bf 134} (1998), no. 2,  301--333. Erratum: Invent. Math.  {\bf 138}  (1999),  no. 1, 225.

\Ref[F]
G. Faltings,
\sl Integral crystalline cohomology over very ramified
valuation rings,
\rm J. Amer. Math. Soc. {\bf 12} (1999), no. 1,  117--144.

\Ref[FC]
G. Faltings and C.-L. Chai,
\sl Degeneration of abelian varieties,
\rm  Ergebnisse der Mathematik und ihrer Grenzgebiete (3), Band {\bf 22}, Springer-Verlag, Berlin, 1990.

\Ref[G1] 
A. Grothendieck, 
\sl \'El\'ements de g\'eom\'etrie alg\'ebrique. IV. \'Etude locale des sch\'emas et des morphismes de sch\'ema (troisi\`eme and quatrie\`eme Parties), 
\rm Inst. Hautes \'Etudes Sci. Publ. Math., Vol. {\bf 28} (1966) and Vol. {\bf 32} (1967).

\Ref[G2]
A. Grothendieck et al.,
\sl Cohomologie locale des faisceau coh\'erents et th\'eor\`ems de Lefschetz locaux et globaux,
\rm S\'eminaire de G\'eom\'etrie Alg\'ebrique du Bois-Marie, 1962. Advanced Studies in Pure Mathematics, Vol. {\bf 2}, North-Holland Publishing Co., Amsterdam; Masson \& Cie, \'Editeur, Paris, 1968.

\Ref[H1]
J. E. Humphreys, 
\sl Introduction to Lie algebras and representation theory,
\rm Grad. Texts in Math., Vol. {\bf 9}, Springer-Verlag, New York-Berlin, 1972.

\Ref[H2]
J. E. Humphreys, 
\sl Conjugacy classes in semisimple algebraic groups, 
\rm Mathematical Surveys and Monographs, Vol. {\bf 43}, Amer. Math. Soc., Providence, RI, 1995.  

\Ref[He]
S. Helgason,
\sl Differential geometry, Lie groups, and symmetric spaces,
\rm  Pure and Applied Mathematics, Vol. {\bf 80}, Academic Press, Inc., New York-London, 1978.

\Ref[J]
J. C. Jantzen,
\sl Representations of algebraic groups. Second edition,
\rm Mathematical Surveys and Monographs, Vol. {\bf 107}, Amer. Math. Soc., Providence, RI, 2003.

\Ref[M] 
H. Matsumura, 
\sl Commutative algebra. Second edition, 
\rm Mathematics Lecture Note Series, Vol. {\bf 56}, The Benjamin/Cummings Publ. Co., Inc., Reading, MA, 1980.

\Ref[MB]
L. Moret-Bailly,
\sl Un th\'eor\`eme de puret\'e pour les familles de courbes lisses,
\rm C. R. Acad. Sci. Paris S\'er. I Math.  {\bf 300}  (1985),  no. 14,  489--492.
\Ref[MS]
W. H. Mills and G. B. Seligman, 
\sl Lie Algebras of classical Lie type,
\rm  J. Math. Mech. {\bf 6} (1957), 519–-548.

\Ref[Mi]
W. H. Mills,
\sl Classical type Lie algebras of characteristic 5 and 7,
\rm J. Math. Mech. 6 (1957), 559-–566.

\Ref[PY]
G. Prasad and J.-K. Yu, 
\sl On quasi-reductive group schemes,
\rm J. Algebraic Geom.  {\bf 15}  (2006),  no. 3,  507--549.

\Ref[S]
G. B. Seligman,
\sl Modular Lie Algebras,
\rm Ergebnisse der Mathematik und ihrer Grenzgebiete, Band {\bf 40}, Springer-Verlag, New-York, 1967.

\Ref[V1]
A. Vasiu,
\sl Integral canonical models of Shimura varieties of preabelian type,
\rm Asian J. Math. {\bf 3} (1999), no. 2,  401--518.
 
\Ref[V2] 
A. Vasiu,
\sl A purity theorem for abelian schemes,
\rm Michigan Math. J. {\bf 52} (2004), no. 1,  71--81.

\Ref[V3] 
A. Vasiu,
\sl On two theorems for flat, affine group schemes over a discrete valuation ring,
\rm Centr. Eur. J. Math. {\bf 3} (2005), no. 1,  14--25. 

\Ref[V4] 
A. Vasiu,
\sl Normal, unipotent subgroup schemes of reductive groups,
\rm C. R. Math. Acad. Sci. Paris {\bf 341} (2005), no. 2,  79--84.

\Ref[V5] 
A. Vasiu,
\sl Integral models in unramified mixed characteristic (0,2) of hermitian orthogonal Shimura varieties of PEL type, Part I,
\rm J. Ramanujan Math. Soc. {\bf 27} (2012), no. 4, 425--477.

\Ref[V6] 
A. Vasiu,
\sl Good reductions of Shimura varieties of Hodge type in arbitrary unramified mixed characteristic, Parts I and II,
\rm math 0707.1668 and math 0712.1572.

\Ref[V7] 
A. Vasiu,
\sl Three methods to prove the existence of integral canonical models of Shimura varieties of Hodge type
\rm work in progress, see also http://arxiv.org/abs/0811.2970

\Ref[VZ]
A. Vasiu and T. Zink, 
\sl Purity results for $p$-divisible groups and abelian schemes over regular bases of mixed characteristic, 
\rm Doc. Math. {\bf 15} (2010), 571--599.

}}

\bigskip
\hbox{Adrian Vasiu,}
\hbox{Dept. of Mathematical Sciences, Binghamton University,}
\hbox{P.O. Box 6000, Binghamton, NY 13902-6000, U.S.A.}
\hbox{e-mail: adrian\@math.binghamton.edu}

\enddocument